\documentclass[format=acmsmall,manuscript,screen]{acmart}
\pdfoutput=1

\newif\ifarxiv
\arxivtrue

\setcopyright{none}

\usepackage[utf8]{inputenc}
\usepackage{array}
\usepackage{algorithm}
\usepackage{algpseudocode}
\usepackage{amsmath}
\usepackage{booktabs}
\usepackage[english]{babel}
\usepackage{enumitem}
\usepackage{hyperref}
\usepackage{pgfplots}
\usepackage{pgfplotstable}
\usepackage{cleveref}
\usepackage{tikz}
\usetikzlibrary{calc}
\usepackage{verbatim}
\usepackage{subcaption}
\usepackage{xspace}

\usetikzlibrary{shapes,arrows,positioning,decorations.pathreplacing,
        intersections,matrix,patterns,fadings,calc,decorations.markings,arrows.meta}

\pgfplotsset{compat=1.17}

\newcommand{\bfc}{\mathbf{c}}
\newcommand{\bfb}{\mathbf{b}}
\def\hdiv{{H(\mathrm{div})}}

\newcommand{\Irksome}{\texttt{Irksome}}

\usepackage{listings}
\definecolor{DarkBlue}{rgb}{0.00,0.00,0.55}
\definecolor{DarkRed}{rgb}{0.55,0.00,0.00}
\definecolor{DarkGreen}{rgb}{0.00,0.55,0.00}
\definecolor{Bittersweet}{rgb}{1.0, 0.44, 0.37}
\definecolor{Purple}{rgb}{0.5, 0.0, 0.5}

\acmJournal{TOMS}

\title{Extending Irksome: improvements in automated Runge--Kutta time stepping for finite element methods}

\author{Robert C. Kirby}
\affiliation{%
  \institution{Baylor University}
  \department{Department of Mathematics}
  \streetaddress{1410 S.~4th St.}
  \city{Waco}
  \state{TX}
  \country{USA}}
\email{robert_kirby@baylor.edu}

\author{Scott P. MacLachlan}
\affiliation{%
  \institution{Memorial University of Newfoundland}
  \department{Department of Mathematics and Statistics}
  \city{St.~John's}
  \state{NL}
  \country{Canada}
}
\email{smaclachlan@mun.ca}

\numberwithin{equation}{section}
\crefname{algorithm}{Algorithm}{Algorithms}
\crefname{figure}{Fig.}{Figs.}
\crefname{table}{Table}{Tables}

\begin{abstract}
  \texttt{Irksome} is a library based on the Unified Form Language (UFL) that enables automated generation of Runge--Kutta methods for time-stepping finite element spatial discretizations of partial differential equations (PDEs).
  Allowing users to express semidiscrete forms of PDEs, it generates UFL representations for the stage-coupled variational problems to be solved at each time step.
  The Firedrake package then generates efficient code for evaluating these variational problems and allows users a wide range of options to deploy efficient algebraic solvers in PETSc.
  In this paper, we describe several recent advances in \texttt{Irksome}.
  These include alternate formulations of the Runge--Kutta time-stepping methods
  and optimized support for diagonally implicit (DIRK) methods.  Additionally, we present new and improved tools for building preconditioners for the resulting linear and linearized systems, demonstrating that these can lead to efficient approaches for solving fully implicit Runge-Kutta discretizations.
  The new features are demonstrated through a sequence of computational examples demonstrating the high-level interface and obtained solver performance.
\end{abstract}


\begin{document}
\maketitle

\section{Introduction}
High-level finite element software like FEniCS~\cite{Logg:2012}, deal.II~\cite{Bangerth:2007}, and Firedrake~\cite{Rathgeber:2016} has long been more successful at providing abstractions for spatial discretizations than temporal ones.
Many codes leave users to manually write loops over time steps or interface to an external time-stepping package.
However, recent literature on high-order approximation of partial differential equations (PDEs) reveals renewed interest in higher-order fully implicit Runge--Kutta (RK) methods, which are difficult to implement and not typically included in
libraries for the solution of ordinary differential equations (ODEs).
Fully implicit methods offer strong stability properties and high accuracy, at the cost of solving a complicated algebraic system coupling the method stages at each time step.
Diagonally implicit methods are cheaper to compute but may suffer from reduced accuracy due to low stage order.
On the other hand, fully implicit methods can have higher stage order and hence better overall accuracy.
Therefore, there is interest in having effective software and solvers to make fully implicit methods more practical.


Realizing the potential of fully implicit methods in practice requires effective solution strategies for the stage-coupled algebraic system to be solved at each time step.
Among the earliest work on this problem is that of~\citet{staff2006preconditioning, mardal2007order, nilssen2011order}, where block-structured preconditioners are developed and analysed, and of~\citet{vanlent2005}, where multigrid methods are applied directly to the stage-coupled system.
The basic motivation for block-structured preconditioners is quite straightforward, even if the analysis is involved -- the Jacobian of the stage-coupled system is block-structured, and one approximates it with its block-diagonal or block-triangular part in the fashion of block Jacobi or Gauss-Seidel.
These methods allow one to use a suitable preconditioner for a backward Euler method (e.g. multigrid) to approximate the inverse of the diagonal blocks; however, regardless of how accurately this is done, the methods tend to require more outer iterations as the stage count increases.
More powerful techniques that allow for stage-independent convergence rates for a wider range of problems~\cite{masud2021new,southworth2022fast1,southworth2022fast2,abu2022monolithic} have also been developed, and we will describe such techniques later in Section~\ref{sec:precond}.

Effective software for applying RK methods to challenging PDEs problems must streamline both the construction and algebraic solution process for fully implicit methods. \Irksome{} attempts to do just this within the Firedrake project~\cite{Rathgeber:2016, kirby2018solver}.
In our previous work~\cite{farrell2021irksome}, we introduced the \Irksome{} package for combining fully implicit Runge-Kutta time stepping with spatial finite element discretizations.
This package is based on UFL (United Form Language)~\cite{Alnaes:2014} manipulation -- after introducing a new node representing time derivatives, we transform the semidiscrete symbolic problem description into UFL for the per-time-step multi-stage coupled system.
Then, \Irksome{} allows these problems to be solved using PETSc~\cite{petsc-user-ref, petsc-efficient, petsc-web-page} and provides rules for updating the solution.
It works on mixed (product of multiple approximating spaces) and nonlinear problems, and also allows users to deploy many effective solver strategies for the stage-coupled system.  In this work, we introduce several important extensions to the \Irksome{} package, namely:
\begin{enumerate}
\item High-accuracy implementation of strong Dirichlet boundary conditions as algebraic equations coupled with the given differential equation;
\item Stage value formulations~\cite{hairer2006geometric} and alternative stage derivative formulations~\cite{butcher1976implementation} to complement the original stage derivative formulation;
\item Optimized support for diagonally implicit Runge-Kutta (DIRK) methods, building and reusing a single stage solver rather than the fully coupled methodology needed for fully implicit RK methods; and
\item Improved support for natural preconditioners for fully implicit RK discretizations.
\end{enumerate}

The remainder of this paper is outlined as follows.  \Cref{sec:RK} reviews the Runge-Kutta methodology as applied to space-time PDEs that are semi-discretized using finite element methods.
In this section, we introduce the alternative stage formulations mentioned above, treatment of boundary conditions as algebraic components, and DIRK schemes.
Tailored preconditioners for RK methods and their realization in \Irksome{} are discussed in~\Cref{sec:precond}.  
Finally, numerical examples (corresponding to new documented example codes in the \Irksome{} library) are presented in~\Cref{sec:examples}, followed by conclusions and outlook in~\Cref{sec:conclusions}.

\section{Runge--Kutta methods}\label{sec:RK}
Runge--Kutta methods were originally formulated for ordinary  differential equations of the form
\begin{equation}
  \label{eq:ode}
y^\prime(t) - F(t, y) = 0,
\end{equation}
where $F: (0,T] \times \mathbb{R}^{m} \rightarrow \mathbb{R}^m$,
and the solution $y: (0,T] \rightarrow \mathbb{R}^m$. The equation must also
satisfy some initial condition
\begin{equation}
  y(0) = y^0.
\end{equation}
Given some approximation to the solution $y^n \approx y(t^n)$ and some
$t^{n+1} = t^n + \Delta t$,
Runge--Kutta methods approximate $y(t^{n+1})$ by
\begin{equation}
  y^{n+1} = y^n + \Delta t \sum_{i=1}^s b_i k_i,
\end{equation}
where, for all $1 \leq i \leq s$, the \emph{stages} $k_i \in \mathbb{R}^m$ satisfy
\begin{equation}
  \label{eq:stageeq}
  k_i - F\left(t^n+c_i\Delta t, y^n+\Delta t\sum_{j=1}^s a_{ij} k_j
  \right) = 0.
\end{equation}
The numbers $a_{ij}$, $b_i$, and $c_i$ for $1 \leq i, j \leq s$ are, in
principle, arbitrary, but are chosen so that the resulting method has a
given order of accuracy as well other desired properties (e.g.~various
notions of stability or symplecticity).  They are typically organized
in a \emph{Butcher tableau}
\begin{equation}
  \begin{array}{c|c}
    \bfc & A \\ \hline
    & \bfb
  \end{array},
\end{equation}
where the vectors $\bfb,\bfc \in \mathbb{R}^s$ contain the entries $b_i$ and
$c_i$, respectively, and $A \in \mathbb{R}^{s \times s}$ contains the
entries $a_{ij}$.  If $A$ is strictly lower triangular, then the method is explicit -- each stage value can be computed in sequence without recourse to an algebraic system (modulo mass matrices in the variational context). Otherwise, the method is implicit.  In the case of a fully implicit method ($A$ being essentially dense), one must solve an $(ms) \times (ms)$ system of algebraic equations to determine all the stage values simultaneously.  When $A$ is lower triangular but not strictly so (a diagonally implicit method or DIRK), one may solve $s$ consecutive $m \times m$ algebraic systems for the stages.
This has historically been regarded as more efficient than solving the single large system for fully implicit methods, but may come at a cost of reduced accuracy in the approximation of $y^{n+1}$ compared to fully implicit schemes with the same number of stages.
Additionally, recent work~\cite{farrell2021irksome, masud2021new,southworth2022fast1,southworth2022fast2,abu2022monolithic}, as well as the examples in this paper, suggest that fully implicit methods may be competitive in run-time, especially relative to the accuracy obtained.

We note two important aspects of the RK methodology not mentioned above.
First of all, when Equation~\eqref{eq:ode} comes from applying the method-of-lines methodology to a spatial finite element discretization, the most natural formulation has the finite element mass matrix multiplying the time-derivative term in this equation.
To write the system in this form requires pre-multiplying through by the inverse of that mass matrix, which is computationally feasible using sparse direct or iterative solvers, but not always a trivial computational task.
Since our primary interest is in implicit methods, where the mass matrix is just one term in the (non-)linear systems to be solved at each time step, we do not do this pre-multiplication explicitly.  Rather, we consider a generalized form for Equation~\eqref{eq:ode}.
Secondly, we note that RK methods can also be formulated in a similar way for differential-algebraic equations (DAEs), as was considered in~\cite{wanner1996solving, abu2022monolithic}.
We expand on this below.

To fix ideas, we will consider two examples in some detail.  We first apply a general RK method to a finite element discretization of the heat equation, finding $u(t,\cdot)\in V_h$ such that 
\begin{equation}
\label{eq:heat}
  \left(u_t, v\right) - \left(\nabla u, \nabla v\right) = \left(f, v \right),
\end{equation}
for all $v\in V_h$, posed on some domain $\Omega \subset \mathbb{R}^d$ with $d \in \{1,2,3\}$,
together with Dirichlet boundary conditions,
$u|_{\partial \Omega} = g(t, \cdot)$.  We let $V_h$ consist of standard continuous
piecewise polynomials defined over a triangulation of $\Omega$.
Applying a generic $s$-stage RK method to~\eqref{eq:heat} leads to a variational problem for the $s$ stages.  We seek $\{ k_i \}_{i=1}^s \subset V_h$ such that
\begin{equation}
\label{eq:rkheat}
\left( k_i , v_i \right) + \left( \nabla \left( u^n + \Delta t \sum_{j=1}^s a_{ij} k_j \right), \nabla v_i \right) 
- \left(f\left(t^n + c_i \Delta t , \cdot \right), v_i \right)= 0,
\end{equation}
for all $v_i \in V_h$, $1 \leq i \leq s$, and then define $u^{n+1}$ as
\begin{equation}
u^{n+1} = u^n + \sum_{i=1}^s b_i k_i.
\end{equation}

As an example, we consider the two-stage LobattoIIIC method, given by Butcher tableau
\begin{equation}
  \label{eq:lobattotableau}
  \begin{array}{c|cc}
    0 & {1}/{2} & -{1}/{2}  \\
    1 & {1}/{2} & {1}/{2} \\ \hline
    & {1}/{2} & {1}/{2}
  \end{array},
\end{equation}
which leads to the variational problem on $V_h \times V_h$
\begin{equation}
  \label{eq:heatvarrk}
  \begin{split}
    \left( k_1 , v_1 \right)
    + \left( \nabla \left( u^n + \Delta t \left( \tfrac{1}{2} k_1 - \tfrac{1}{2} k_2 \right) \right), \nabla v_1 \right) 
    - \left(f\left(t , \cdot \right), v_1 \right) & = 0, \ \ \ \forall v_1 \in V_h, \\        \left( k_2 , v_2 \right) + \left( \nabla \left( u^n + \Delta t \left( \tfrac{1}{2} k_1 + \tfrac{1}{2} k_2 \right) \right), \nabla v_2 \right) 
    - \left(f\left(t + \Delta t , \cdot \right), v_2 \right) & = 0, \ \ \ \forall v_2 \in V_h,
  \end{split}
\end{equation}
with $u^{n+1} = u^n + \frac{\Delta t}{2}(k_1+k_2)$.  The critical contribution of \Irksome{} is that it takes the Butcher tableau and a UFL description of the semidiscrete problem in~\eqref{eq:heat} and generates a UFL description of~\eqref{eq:heatvarrk}.
This is done by manipulating the expression graph, as described in our earlier work~\cite{farrell2021irksome}.

Posing Dirichlet boundary conditions on the $k_i$ such that the computed solution, $u^{n+1}$, will satisfy $u^{n+1}(t,\cdot) = g(t^{n+1}, \cdot)$ on the boundary to sufficient accuracy is a challenge.
Previously, in~\cite{farrell2021irksome}, we proposed an heuristic based on the notion that the $k_i$ approximate time derivatives at $t^n + c_i \Delta t$, suggesting that $k_i|_{\partial\Omega} = g'(t^{n+1}+c_i\Delta t,\cdot)$.
However, this approach has proven to have some inherent weaknesses, and similar approaches are well-known to limit the overall accuracy of the time integrator~\cite{doi:10.1137/0916072}.  As described below, we now adopt a DAE-viewpoint on enforcing Dirichlet boundary conditions that retains (at least) the stage order of the scheme, which matches the expected convergence of the time-stepper when applied to stiff ODEs or to DAEs.

Many important PDEs lead to differential algebraic rather than ordinary differential equations after spatial discretization.
For our second example, we consider such a system in the incompressible Navier-Stokes equations,
\begin{equation}
  \label{eq:nse}
  \begin{split}
    \mathbf{u}_t + \mathbf{u} \cdot \nabla \mathbf{u}  - \tfrac{1}{Re} \Delta \mathbf{u} + \nabla p & = \mathbf{0}, \\
    \nabla \cdot \mathbf{u} & = 0,
  \end{split}
\end{equation}
where $\mathbf{u}$ is the fluid velocity, $p$ is the pressure, and the dimensionless Reynolds number, $Re$, measures the ratio of inertial to viscous forces.  Typically,
these equations are posed in some bounded domain in $\mathbb{R}^d$ with $d=2, 3$, and the system is driven by boundary conditions,
although driving forces (e.g. due to gravity) in the first equation are not uncommon.
Such a system requires initial data for the velocity, $\mathbf{u}(0, \mathbf{x}) = \mathbf{g}(\mathbf{x})$ for $\mathbf{x}\in\Omega$.  When integrating with an RK timestepper, we also require initial data for the pressure, $p(0, \mathbf{x}) = q(\mathbf{x})$ for $\mathbf{x}\in\Omega$, as described below.
By choosing appropriate finite element spaces, $V_h$ and $W_h$, for velocity and pressure, a standard Galerkin method for~\eqref{eq:nse} is to seek
$\mathbf{u} \in V_h$ and $p \in W_h$ such that
\begin{equation}
  \label{eq:nsevarform}
  \begin{split}
    \left( \mathbf{u}_t, \mathbf{v} \right) \
    + \left( \mathbf{u} \cdot \nabla \mathbf{u}, \mathbf{v} \right)
    + \tfrac{1}{Re} \left( \nabla \mathbf{u}, \nabla \mathbf{v} \right)
    - \left( p, \nabla \cdot \mathbf{v} \right) & = 0, \ \ \ \forall \mathbf{v} \in V_h \\
    \left( \nabla \cdot \mathbf{u}, w \right) & = 0, \ \ \ \forall w \in W_h.
  \end{split}
\end{equation}
Readers familiar with the Navier-Stokes equations will note that the particular choices of $V_h$ and $W_h$ have a significant impact on stability, accuracy, and efficiency of the finite-element discretization.
While the same concerns arise for the time discretization, the choice of a suitable time integrator is largely independent of that of the spatial discretization.

We also note that the assumption of Galerkin discretization with inf-sup approximating spaces is primarily illustrative -- other effective spatial discretizations (e.g. stabilized or discontinuous Galerkin) may also be advanced in time using the approaches we describe here.

The two-stage LobattoIIIC method given by~\eqref{eq:lobattotableau} can be applied to~\eqref{eq:nsevarform}, giving a variational problem on
$Z_h = V_h \times W_h \times V_h \times W_h$.
Here, we have unknowns for the velocity and pressure for each stage.
We seek $\left(\mathbf{k}_1^\mathbf{u}, k_1^{p}, \mathbf{k}_2^\mathbf{u}, k_2^{p}\right) \in Z_h$ such that, with
\begin{equation}
    \mathbf{u}_i = \mathbf{u}^n + \Delta t \sum_{j=1}^s a_{ij} \mathbf{k}^\mathbf{u}_j, \quad
    p_i = p^n + \Delta t \sum_{j=1}^s a_{ij} k^p_j,
\end{equation}
we solve
\begin{equation}
  \label{eq:nsevarformrk}
  \begin{split}
    \left( \mathbf{k}^\mathbf{u}_1, \mathbf{v}_1 \right) \
    + \left( \mathbf{u}_1 \cdot \nabla \mathbf{u}_1, \mathbf{v}_1 \right)
    + \tfrac{1}{Re} \left( \nabla \mathbf{u}_1, \nabla \mathbf{v}_1 \right)
    - \left( p_1, \nabla \cdot \mathbf{v}_1 \right) & = 0, \ \ \ \forall \mathbf{v}_1 \in V_h \\
    \left( \nabla \cdot \mathbf{u}_1, w_1\right) & = 0, \ \ \ \forall w_1 \in W_h, \\
    \left( \mathbf{k}^\mathbf{u}_2, \mathbf{v}_2 \right) \
    + \left( \mathbf{u}_2 \cdot \nabla \mathbf{u}_2, \mathbf{v}_2 \right)
    + \tfrac{1}{Re} \left( \nabla \mathbf{u}_2, \nabla \mathbf{v}_2 \right)
    - \left( p_2, \nabla \cdot \mathbf{v}_2 \right) & = 0, \ \ \ \forall \mathbf{v}_2 \in V_h \\
    \left( \nabla \cdot \mathbf{u}_2, w_2\right) & = 0, \ \ \ \forall w_2 \in W_h.
  \end{split}
\end{equation}
Here, the unknowns are $\mathbf{k}_1^\mathbf{u}, k_1^{p}, \mathbf{k}_2^\mathbf{u},$ and $k_2^{p}$, and the values $\mathbf{u}_i$ and $p_i$ are defined algebraically from these as above.
The variational problem~\eqref{eq:nsevarformrk} couples all of the stages, since each of $\mathbf{u}_i$, $p_i$ contains the unknowns for both stages.
Again, \Irksome{} maps a UFL description of~\eqref{eq:nsevarform} into a UFL description for~\eqref{eq:nsevarformrk}. 
We note that when considering accuracy of the integrator for degree-2 DAEs, such as the semi-discrete equations in~\eqref{eq:nsevarform}, it is the so-called \emph{stage order} of the scheme that matters, and not the global order of the scheme~\cite{wanner1996solving}.

We give an example of \Irksome{}'s interface in Figure~\ref{fig:nsecode}, where we pose a Taylor-Hood discretization of Navier-Stokes for the lid-driven cavity problem.
Besides the UFL extension for time derivatives \lstinline{Dt} and classes for Butcher tableaux like \lstinline{RadauIIA}, the main entry point for users is \lstinline{TimeStepper}.
Originally, this was a class that constructed the variational problem to solve for each time step, configured a PETSc solver, and provided methods to advance the solution forward in time.
In order to accommodate the new features described in this work with a backward-compatible interface, we have evolved \lstinline{TimeStepper} to be a factory function with keyword arguments used to select or configure the new features.  We note that the code in Figure~\ref{fig:nsecode} is intentionally simple, and that the example codes for this paper~\cite{zenodo/Firedrake-20250218.0} include more details on constructions appropriate for running in parallel with performant solver parameters.

\begin{figure}
  \begin{lstlisting}
from firedrake import *
from Irksome import Dt, TimeStepper, RadauIIA

msh = UnitSquareMesh(16, 16)
V = VectorFunctionSpace(msh, 'CG', 2)
W = FunctionSpace(msh, 'CG', 1)
Z = V * W
up = Function(Z)
u, p = split(up)
v, w = TestFunctions(Z)

Re = Constant(10.0)

t = Constant(0.0)
dt = Constant(1.0 / 16)

F = (inner(Dt(u), v) * dx + inner(dot(u, grad(u)), v) * dx
     + 1/Re * inner(grad(u), grad(v)) * dx - inner(p, div(v)) * dx
     + inner(div(u), w) * dx)

bcs = [DirichletBC(Z.sub(0), 0, (1, 2, 3)),
       DirichletBC(Z.sub(1), as_vector([1, 0]), (4,))]
     
butcher_tableau = RadauIIA(2)
stepper = TimeStepper(F, butcher_tableau, t, dt, up, bcs=bcs)
  \end{lstlisting}
  \caption{Basic interface for implementing Navier-Stokes in \Irksome{}.}
  \label{fig:nsecode}
\end{figure}

\subsection{Treatment of Dirichlet boundary conditions}

The DAE perspective discussed above for the Navier-Stokes equations provides a methodology for ensuring accurate integration of strongly enforced Dirichlet boundary conditions.  Given a finite element space $V_h$, we consider a standard decomposition of a function $u\in V_h$ into the basis functions of $V_h$ that are non-zero on the boundary in question (either all of $\partial\Omega$ or a specified subset, $\Gamma \subset\partial\Omega$) and those that are not.  Posing a boundary condition on $u(t,\cdot)$ on segment $\Gamma$ is equivalent to posing an algebraic equation for the coefficients of the basis functions associated with that segment.

Consider the Dirichlet BC that $u(t,\cdot) = g(t,\cdot)$ on $\Gamma\subset\partial\Omega$.  If we restrict $u$ to $\Gamma$, we can view the boundary condition as requiring that
\begin{equation}
u_i \approx u(t^n + c_i \Delta t,\cdot)  = g(t^n + c_i \Delta t,\cdot) \text{ on }\Gamma,
\end{equation}
for $1\leq i \leq s$, where $u_i = u^n + \Delta t\sum_{i=1}^s a_{ij}k_j$ is the approximation of $u(t,\cdot)$ at the stage time $t^n + c_i \Delta t$.  Setting $u_i$ equal to the boundary value, we get a system of equations, with
\begin{equation}
u^n + \Delta t \sum_{j=1}^s a_{ij} k_j = g(t^n + c_i \Delta t,\cdot) \text{ on }\Gamma,
\end{equation}
for $1 \leq i \leq s$.  We rewrite this as
\begin{equation}
\sum_{j=1}^s a_{ij} k_j = \frac{1}{\Delta t} \left(g(t^n + c_i \Delta t,\cdot) - u^n\right)  \text{ on }\Gamma,
\end{equation}
for $1 \leq i \leq s$.  When the Butcher matrix, $A$, is invertible, we can solve this system of equations exactly for the boundary conditions on $k_i$, at least up to error in the $L^2$ projection of $g(t^n + c_i \Delta t,\cdot)$ into $V_h$.  For the fully implicit schemes considered here, this condition is easily seen to be satisfied, and we refer to implementing boundary conditions of this form as using ``DAE-type boundary conditions''.  We note that the boundary values prescribed on the stages do not guarantee that $u^{n+1}$ exactly equals $g(t^{n+1},\cdot)$ on $\Gamma$, since the conditions are imposed on the stage approximations and not the reconstructed solution, $u^{n+1}$.  For stiffly accurate schemes, however, including the Lobatto IIIC and Radau IIA schemes considered here, since $u^{n+1} = u_s$, we get the additional guarantee that $u^{n+1}$ does match its prescribed Dirichlet boundary data (up to projection error).

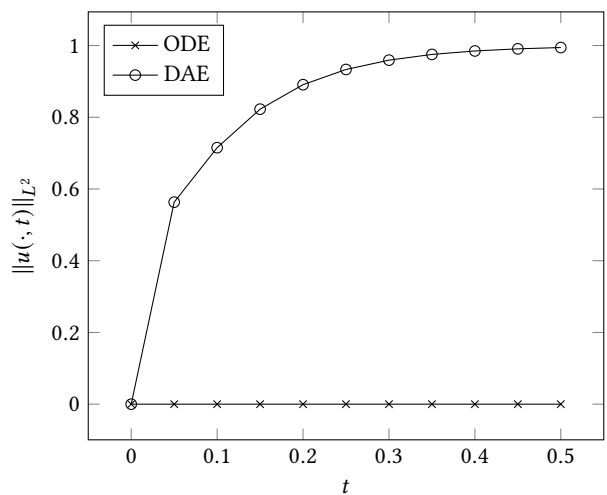
\begin{figure}
  \begin{center}
    \begin{tikzpicture}
      \begin{axis}[xlabel= $t$, ylabel={$\|u(\cdot, t)\|_{L^2}$},
        legend pos=north west]
        \addplot[mark=x]
        table[x=t,y=nrmu,col sep=comma]{odenorm.csv};
        \addlegendentry{ODE};
        \addplot[mark=o]
        table[x=t,y=nrmu,col sep=comma]{daenorm.csv};
        \addlegendentry{DAE};        
      \end{axis}
    \end{tikzpicture}
  \end{center}
  \caption{Comparing ODE and DAE boundary condition enforcement for the heat equation with incompatible initial and boundary data.}
  \label{fig:heatbc}
\end{figure}

While this implementation of Dirichlet boundary conditions is more consistent with Runge-Kutta formulations, we also note that incompatible boundary and initial data can lead to difficulties for the original ``ODE-type'' formulation described in~\cite{farrell2021irksome}.
For example, consider the one-dimensional heat equation $u_t - u_{xx} = 0$ on $[0, 1]$ with initial condition $u(0, x) = 0$ and Dirichlet boundary conditions $u(t, 0) = u(t, 1) = 1$ for $t > 0$.
A simple calculation (for example, by separation of variables) shows that the exact solution rapidly converges to $u(x, t) = 1$ as $t$ grows.
Discretizing this problem with standard Lagrange finite elements and an implicit time-stepper captures this behavior using DAE-type boundary conditions.
However, because the ODE-type boundary conditions consider the time derivative of the boundary data (which is zero in this case), they miss the incompatibility and produce a steady zero solution.
To illustrate this, we took the finite element space comprising piecewise linears over ten intervals and integrated the heat equation until time $0.5$ using the three-stage LobattoIIIC method with time steps of size $0.05$.
Figure~\ref{fig:heatbc} plots the $L^2$ norm of the finite element solution at each time using both ODE- and DAE-type boundary conditions.

Both approaches remain available in \Irksome{}.
One specifies the method of enforcing boundary conditions with the optional keyword argument \lstinline{bc_method} of \lstinline{TimeStepper}.
It takes either the string \lstinline{'DAE'} (the method described here) or \lstinline{'ODE'} (the method originally described in~\cite{farrell2021irksome}).
If no keyword argument is specified, the default of \lstinline{'DAE'} is taken.  We note that many practical cases show little difference between the two approaches, but the DAE approach is both more robust and supported by standard RK convergence theory.

We note that another alternative for implementing Dirichlet boundary conditions is through weak enforcement, using penalty methods, such as Nitsche's method~\cite{Nitsche_1971, benzaken2022constructing}.  While this avoids the potential difficulties discussed above, it brings its own complications, particularly in determining a penalty formulation that does not adversely affect the error estimates of the spatial discretization.  These methods are incorporated in the variational form \lstinline{F} and not through a \lstinline{firedrake.DirichletBC} object or the choice of the keyword argument.

\subsection{Stage formulations}

Our prior work dealt with the standard formulation of RK methods, which we refer to as the \emph{stage derivative} formulation, but other formulations have long been known in the literature.
We return to the generic ODE~\eqref{eq:ode} to describe these approaches, although they are readily adapted to variational problems and DAEs.
Going back to~\citet{butcher1976implementation}, a change of variables is possible when the matrix $A$ is invertible, solving for linear combinations of the $k_i$ which can then be separated after the system is solved.
This leads to different algebraic structure in the system to be solved at each time step, which can have algorithmic advantages.

To introduce these schemes, we define new variables $w_i \in V_h$ by
\begin{equation}
  \label{eq:butcherw}
  w_i = \sum_{j=1}^s a_{ij} k_j,
\end{equation}
which allows us to rewrite~\eqref{eq:stageeq} as
\begin{equation}
  \label{eq:stageeqw}
  \sum_{j=1}^s \left(\left(A^{-1}\right)_{ij} w_j\right) - F\left(t^n+c_i\Delta t, y^n+\Delta t w_i
  \right) = 0,
\end{equation}
for $1\leq i \leq s$.
In this form, the coupling between stage values appears in the first term, outside of the nonlinearity in $F$; the contribution of the nonlinearity to the Jacobian is logically block diagonal.
This can significantly reduce the cost of assembling and possibly solving the Jacobian since there are $s$ rather than $s^2$ blocks of the matrix in which it appears.
However, if the overall process is dominated by solution rather than assembly, this may not be significant.
Our initial implementation of \Irksome{} described in~\cite{farrell2021irksome} focused on generality - working for a wide range of RK methods.
However, we have updated our implementation to allow the user to provide an option to switch between~\eqref{eq:stageeq} and~\eqref{eq:stageeqw}, with an internal check for the invertibility of $A$.
Further generalizations are also possible -- we may write $A = A_1A_2$ for any invertible matrix $A_2$, then define the stage values $w = A_2 k$, although we are unaware of other splittings that provide similar benefits as taking $A_2 = A$ with $A_1 = I$.

Like methods of boundary condition enforcement, the choice of splitting is specified by a keyword argument \lstinline{splitting} to \lstinline{TimeStepper}, which can be either \lstinline{irksome.tools.IA} or \lstinline{irksome.tools.AI}.  The default option is \lstinline{AI}, leading to the classical formulation in~\eqref{eq:stageeq}.
No other user modification is required to switch between these choices.

Another alternative is to formulate RK methods in terms of the unknown values of the solution at various times in each interval.
Returning to~\eqref{eq:stageeq}, we could introduce stage values $\{ Y_i \}_{i=1}^s$ satisfying equations
\begin{equation}
  \label{eq:stagevalueform}
  Y_i = y^n - \Delta t \sum_{j=1}^s a_{ij} F(t^n+c_j \Delta t, Y_j),
\end{equation}
and the solution is updated by a linear combination
\begin{equation}
  \label{eq:stagevalueupdate}
  y^{n+1} = y^n - \Delta t \sum_{i=1}^s b_i F(t^n+c_i \Delta t, Y_i).
\end{equation}
While we have written this equivalent formulation for standard ODEs,
we note that, in a finite element context,~\eqref{eq:stagevalueupdate} appears in variational form and requires the inversion of a single-stage mass matrix.
This is a relatively trivial additional cost once the implicit stage-coupled equations have been solved, but must be accounted for in implementation.
For \emph{stiffly accurate} methods, however,~\eqref{eq:stagevalueupdate} simplifies to $y^{n+1}$ taking the final stage value, $Y_s$, which obviates the need for a mass matrix inversion.
We have updated our implementation in \Irksome{} to include this stage value formulation, including automatically detecting and skipping the mass matrix solution for a stiffly accurate method.

Like boundary condition enforcement and splitting strategies, users can switch between stage derivative and stage value formulations by a keyword argument
\lstinline{stage_type} to \lstinline{TimeStepper}, with valid choices being either \lstinline{'deriv'} (the default, leading to~\eqref{eq:stageeq}) or \lstinline{'value'}.
When \lstinline{stage_type='value'} is selected and the Butcher tableau does not represent a stiffly accurate method, the user may also provide
\lstinline{update_solver_parameters} to configure the mass matrix solve required for updating the solution.

\subsection{DIRK schemes}
Although a major motivation in our work has been to make fully implicit RK methods practical, we have also extended \Irksome{} to efficiently implement \emph{diagonally implicit} Runge-Kutta methods.
These methods can provide stability while requiring the solution of a sequence of algebraic problems for each stage rather than that of a single, much larger, stage-coupled system.
However, they historically suffer limitations on accuracy due to their inherent low stage order.
Consequently, one might expect that, if the nontrivial hurdle of solving the algebraic system is overcome, fully implicit methods should provide superior solutions, weighing cost vs. accuracy over a fixed time interval.
Recent work in~\cite{biswas2023design,ketcheson2020dirk} has developed the concept of \emph{weak} stage order, and DIRK schemes with high order and weak stage order are seen to deliver high accuracy in some cases. 

Because of ongoing interest in and current development of new DIRKs, it is important to enable them in our software.
For linear problems, no special treatment is actually required within \Irksome{} to achieve much of the benefit of DIRKs.
Since the stage-coupled algebraic system is linear and block-triangular, it can be solved efficiently by a special choice of PETSc solver options that requires only a single solution with each diagonal block.
Using techniques introduced in~\cite{kirby2018solver}, the large system matrix can be formed in a matrix-free format, only assembling the diagonal blocks (if the choice of solver/preconditioner so requires).
Also, a specialized treatment of DIRKs, optimized to only form and solve the necessary single-stage problems, is now implemented in \Irksome{}.
Moreover, for nonlinear problems, some additional advantage may be gained by solving a sequence of decoupled nonlinear problems for each stage.
These methods are also obtained by selecting an appropriate Butcher tableau with lower-triangular $A$ and passing the keyword argument \lstinline{stage_type='dirk'} to \lstinline{TimeStepper}.
In this case, a single-stage variational problem is created internally and updated between stages.

We have included several DIRKs in \Irksome{} that we will use here.
Among these are the classical method of~\citet{alexander1977diagonally}.
This $L$-stable method has three stages and third-order accuracy, but only has stage order and weak stage order one.
For comparison, we have also included several tableaux from~\cite{biswas2023design, ketcheson2020dirk} with higher weak stage order.
For simplicity, we refer to them by WSODIRK followed by the number of stages, the formal order, and the weak stage order so that WSODIRK433 has four stages, formal order and weak stage order of three.

\section{Preconditioners}\label{sec:precond}
The stage-coupled algebraic systems arising from fully implicit methods require effective solution techniques, and we have focused significant effort on enabling state-of-the-art preconditioners for the resulting block-structured linear systems.
Here, we summarize many approaches known in the literature and describe how to access them in our code stack.

It is well-known that Runge-Kutta time discretization for simple linear PDEs leads to block-structured linear systems of the form
\begin{equation}
  \left( I \otimes M \right) \mathbf{k} + \Delta t \left( A \otimes K \right) \mathbf{k} = \mathbf{f},\label{eq:block_struct_AI}
\end{equation}
where, for the heat equation, $M$ and $K$ are the standard mass and stiffness matrices appearing in a single-stage method.  
Butcher's technique of pulling out~$A^{-1}$ leads to the equivalent linear system
\begin{equation}
  \left( A^{-1}  \otimes M \right) \mathbf{w} + \Delta t \left( I \otimes K \right) \mathbf{w} = \mathbf{f},\label{eq:block_struct_IA}
\end{equation}
with a suitable redefinition of $\mathbf{f}$.

This structure covers most linear evolution-type equations leading to ODEs in semi-discretized form, although DAE-type systems arising from problems like time-dependent Stokes flow require modification as discussed above.
Nonlinear evolution equations require an additional modification of either of these forms.
The stiffness term involving $K$ arises from a Jacobian, so rather than $I \otimes K$ in~\eqref{eq:block_struct_IA}, we have the block diagonal matrix, $\operatorname{diag}(K_1, K_2, \dots, K_s)$, where $K_i$ is a Jacobian obtained by linearizing about stage $i$.
Since the approaches we describe here and our implementation of them in \Irksome{} are indifferent to this issue, we will, for ease of exposition, restrict ourselves to the linear case.
Also, we note that our implementations work seamlessly on the stage value or stage derivative formulations with various splitting strategies

\subsection{Block preconditioners}
Given the blockwise structure of the stage-coupled system, it makes sense to consider block preconditioners.
Block preconditioners have been developed for many kinds of coupled systems, although we focus on approaches specifically applied for Runge-Kutta systems.
Here, we summarize some of these approaches and discuss their implementation in \Irksome{}.
We note that block preconditioners can be readily applied to either of the forms in~\eqref{eq:block_struct_AI} and~\eqref{eq:block_struct_IA}, and present preconditioners for the latter, noting construction of those for the former follows similarly.

Many techniques start from replacing the Butcher matrix $A$ in~\eqref{eq:block_struct_AI} with some $\tilde{A}$ such that the new system is easier to solve (typically block triangular) yet preconditions the original system well.
In some sense, this is a kind of \textit{idealized} preconditioner, in that it still requires exact inversion of the diagonal blocks of the preconditioning matrix.
The techniques in papers such as~\cite{staff2006preconditioning, mardal2007order}, amount to an \emph{additive} approximation of the Butcher matrix $A$.
Writing $A = L + D + U$, in terms of the strictly lower/upper triangular parts, $L$ and $U$, and its diagonal, $D$,
a block diagonal preconditioner for the stage-coupled system in~\eqref{eq:block_struct_IA} is just
\begin{equation}
  \label{eq:pcbjac}
P_{D} = (D^{-1} \otimes M) + \Delta t \left( I \otimes K \right),
\end{equation}
and a block lower-triangular preconditioner would be
\begin{equation}
  \label{eq:pcbgs}
P_{L} = \left( (L + D)^{-1} \otimes M \right) + \Delta t \left( I \otimes K \right).
\end{equation}
Such techniques are block Jacobi and block Gauss-Seidel techniques, respectively, and are attractive because of their simplicity -- they only require a technique to invert the diagonal blocks.
We note that the use of a block lower-triangular preconditioner like $P_L$ is natural for left-preconditioning GMRES; for right-preconditioning either classical GMRES or FGMRES, an upper-triangular preconditioner, $P_{U}$, is naturally defined, replacing $(L+D)^{-1}$ in the first term with $(U + D)^{-1}$.
The empirical results in~\cite{staff2006preconditioning} indicate that right-preconditioning is typically better than left.
Lower-triangular beats block-diagonal preconditioning, which in turn beats upper-triangular.
These results are only for a single model problem.  Other options may be preferred for other problems, but our mathematical and software formalisms support all options with comparable ease.

Under any of these variants, we must provide a scheme for inverting (or at least preconditioning) the diagonal blocks, which take the form $a_{jj}^{-1}M + \Delta t K$.
For a simple problem like the heat equation, where $M$ is the finite element mass matrix and $K$ is the Laplacian stiffness matrix, a natural approach is to approximate the inverse of this block with a standard multigrid method, either geometric or algebraic.
However, even if the blocks are inverted exactly, the number of outer GMRES iterations required to obtain convergence with these preconditioners is observed to grow with the number of stages.

An alternate approach to defining $\tilde{A}$ was introduced in~\citet{masud2021new} and further studied in~\cite{clines2022efficient}.
This turns out to give a far more effective family of block preconditioners.
Rather than the additive splitting $A = L + D + U$ to the Butcher matrix, we begin with a multiplicative decomposition.
We write $A=LDU$  with unit lower triangular $L$, unit upper triangular $U$, and diagonal $D$.
One choice of approximation is to take the lower triangular $\tilde{A}_{LD} = LD$.
Replacing $A$ with the triangular matrix $\tilde{A}_{LD}$ in the coupled system gives a preconditioning matrix for~\eqref{eq:block_struct_IA} as
\begin{equation}
  \label{eq:pcranald}
\tilde{P}_{LD} = \left(\tilde{A}_{LD}^{-1} \otimes M \right) + \Delta t \left( I \otimes K \right).
\end{equation}
Now, $\tilde{P}_{LD}$ is block lower-triangular and applying its inverse only requires inversion of the diagonal blocks, just as with the block Jacobi and Gauss-Seidel approaches above.
However, 
this approach seems to give very low iteration counts, empirically independent of the number of Runge--Kutta stages.

The block Jacobi and Gauss-Seidel methods are readily obtained through the PETSc \lstinline{fieldsplit} mechanism~\cite{brown2012composable}, and we applied them in our earlier paper~\cite{farrell2021irksome}.
However, the Rana-type preconditioners require some degree of care.
Using techniques introduced in~\cite{kirby2018solver}, we give a general implementation of these preconditioners in a programmatic way.
Our preconditioners derive from the abstract class~\lstinline{firedrake.AuxiliaryOperatorPC}, which requires a \lstinline{form} method, providing UFL for a bilinear form to be assembled for a preconditioner.
\Irksome{} places the original semidiscrete UFL form into a PETSc application context, and our preconditioners override the \lstinline{form} method by carrying out the same form manipulations used to derive the algebraic system, but with the Butcher matrix replaced by its approximation $\tilde{A}$.
(Note that this approach makes our implementation quite general, applicable to any problem in \Irksome{}, including nonlinear ones).
This gives a preconditioner of type \lstinline{python} which then allows the user to further configure (approximate) inversion of the resulting matrix, through \lstinline{fieldsplit} or some other means through additional PETSc options.  For example, one can obtain a basic implementation
\begin{lstlisting}
solver_parameters={
  'ksp_type': 'gmres',
  'pc_type': 'python',
  'pc_python_type': 'irksome.RanaLD',
  'aux': {
    'pc_type': 'fieldsplit',
    'pc_fieldsplit_type': 'multiplicative'
  }
}
\end{lstlisting}
This uses default options for inverting the diagonal blocks, and this can be further refined, say, to use algebraic multigrid, as needed.

Other block-structured techniques are also known.
\citet{southworth2022fast1,southworth2022fast2} propose a similar family of preconditioners for~\eqref{eq:block_struct_IA} based on the eigenvalue decomposition of $A^{-1}$.
In the linear case, they use the adjugate expression for the inverse of a block matrix to express a preconditioner that only needs inversion of the matrix $P_s(\Delta t M^{-1}K)$, where $P_s(z)$ is the characteristic polynomial of $A^{-1}$.
Since $A^{-1}$ is $s\times s$, it is not unreasonable to compute its eigenvalue decomposition to high accuracy, and the inverse of $P_s(\Delta t M^{-1}K)$ can be computed by writing it in factored form and inverting the factors termwise.  To avoid complex arithmetic, they propose to treat terms arising from complex-conjugate eigenvalues together, and derive optimal approximate inverses for these blocks.  Each block is, then, inverted using similar multigrid algorithms as discussed above.  While the linear case discussed in~\cite{southworth2022fast1} could be easily implemented in \Irksome{}, the nonlinear case in~\cite{southworth2022fast1} is much more involved, and does not easily map onto PETSc's existing nonlinear solvers.  \citet{leveque2023parallelintime} propose related methods based on the singular value decomposition of $A$, but these are problem-specific and, so far, limited to the linear case.
For this reason, we do not yet support these preconditioners in \Irksome{}.

\subsection{Monolithic multigrid}
The approaches considered so far all hinge on some approximation that decouples the stages, allowing the re-use of single-stage solvers.
The original \Irksome{} paper~\cite{farrell2021irksome} considered an alternative approach of using multigrid schemes that embrace the coupling.
First proposed by~\citet{vanlent2005} for finite-difference discretizations, these \emph{monolithic} multigrid (MMG) approaches rely on relaxation schemes that couple all of the Runge-Kutta stages.
These approaches were further extended to incompressible fluid flow problems in~\cite{abu2022monolithic}.
For a convergence theory of these methods, at least for linear problems, we refer the reader to the recent manuscript~\cite{mmg}.

One may arrive at many useful monolithic methods via an additive Schwarz framework.  Suppose our semidiscrete variational problem is posed on some finite element space, $V$, which we can additively decompose into the sum of $N_p$ spaces by
\begin{equation}
  \label{eq:patches}
V = \sum_{i=1}^{N_p} V_i,
\end{equation}
where the sum need not be direct.  For multigrid relaxation schemes, the spaces $V_i$ will be subspaces of $V$ consisting of functions zero in most of the mesh, perhaps nonzero only in some patch around a vertex.  The natural embedding defines a simple prolongation operator $p_i : V_i \rightarrow V$, and one must also define some restriction $r_i: V \rightarrow V_i$.  Then, an additive Schwarz method specifies the action of an inverse operator $C^{-1}$ by
\begin{equation}
  \label{eq:Cinv}
C^{-1} = \sum_{i=1}^{N_p} p_i A_i^{-1} r_i.
\end{equation}
On each subspace, applying $A_i^{-1} r_i$ amounts to solving a ``local'' problem with a very small number of degrees of freedom -- only the degrees of freedom on a vertex patch, say.

When using lowest-order linear elements for scalar-valued problems on $H^1$, this approach reduces to a point Jacobi relaxation, but it generalizes to higher-order schemes with degree-independent estimates~\cite{Pavarino:1993,Schoeberl:2008} as well, when the subspaces are suitably chosen.
It is also the critical key to obtaining convergent multigrid for problems in $\hdiv$ and $H(\mathrm{curl})$
~\cite{arnold2000multigrid}.

This general framework also helps us to pose monolithic multigrid methods.
If the solution to the PDE lives in some space $V$, the variational problem for the stages will be posed on some $\mathbf{V} = V \times V \times \dots \times V$, and we use the decomposition
\begin{equation}
  \label{eq:bigdecomp}
\mathbf{V} = \sum_{i=1}^{N_p} \mathbf{V}_i,
\end{equation}
where $\mathbf{V}_i = V_i \times V_i \times \dots \times V_i$, and the $V_i$ are the same patches arising in~\eqref{eq:patches}, and one readily defines the analog of~\eqref{eq:Cinv}.  The theory in~\cite{mmg} shows that, provided that multigrid using~\eqref{eq:Cinv} as a relaxation scheme converges for a single stage method (possibly with a complex time step), monolithic multigrid with additive Schwarz relaxation based on the decomposition~\eqref{eq:bigdecomp} also converges.  In fact, the spectral radius of its iteration matrix does not depend explicitly on the number of stages, justifying the largely stage-independent iteration counts we observe for these methods.

Moreover, existing Firedrake features make it simple to deploy such methods.
First of all, the geometric multigrid framework first developed in~\cite{mitchell2016high} has been fully integrated with the PETSc interface described in~\cite{kirby2018solver}.
Second, \lstinline{firedrake.PCPatch}, introduced in~\cite{farrell2021pcpatch}, automates the construction of small-scale, patch-based additive Schwarz methods for UFL-specified variational problems.
A related class, \lstinline{firedrake.ASMPatchPC} performs similar operations in a purely algebraic way, operating on assembled system matrices.
Many options are provided to configure the way in which the various patch-based problems are solved.

\section{Numerical examples}\label{sec:examples}
Our simulations are performed on a 32-core AMD ThreadRipper Pro with 768GB of RAM running Ubuntu 24.04.
We use the configuration of Firedrake, \Irksome{}, PETSc, and other components available through~\cite{zenodo/Firedrake-20240312.0}.
The heat equation results below were all run on 8 cores of the machine,
while the Navier--Stokes and Cahn--Hilliard simulations used 16 cores.

\subsection{The heat equation}
In this section, we compare various fully implicit and DIRK methods for the heat equation.
We pose the problem on the unit square $\Omega = [0, 1]^2$ with exact solution chosen to satisfy $u(t, x, y) = e^{-0.1 t} \sin(\pi x) \cos(\pi y)$, with Dirichlet boundary conditions set to agree with the exact solution on the boundary and imposed using the DAE approach described above.
The domain $\Omega$ is partitioned into an $N \times N$ array of squares, and serendipity elements~\cite{Arnold, crum2022bringing} of degree 2 and 3 are used for the spatial discretization.
With quadratic elements, we expect second and third order convergence in the $H^1$ and $L^2$ norms, respectively, and one order higher in each for cubic elements.

The RadauIIA family provides a theoretically optimal suite of fully implicit methods for this problem in terms of accuracy per stage subject to L-stability.
The $s$-stage method in this family has formal order $2s-1$, and stage order of $s$.
We compare these methods to two DIRK methods.
First, we consider the three-stage, L-stable DIRK in~\cite{alexander1977diagonally}.
While this method has formal third-order accuracy, it has stage order and weak stage order of one and can suffer from order reduction.
We also include the four-stage method WSODIRK433, with formal order and weak stage order of three from~\cite{ketcheson2020dirk}.  Butcher tableaux for these methods are given in~\Cref{sec:butcher}.

As we have discussed above, many options are possible for these formulations.
In the interest of space, we do not pursue an exhaustive analysis here.
To summarize, for a linear, constant-coefficient problem like the heat equation, we found only small differences between classical and stage formulations and IA vs AI splittings, with a slight preference for the IA splitting and stage derivative formulation.
In the sequel, we report performance numbers for this formulation of RadauIIA methods, comparing it to DIRK methods with various solution strategies.

\begin{figure}
  \begin{center}
    \begin{subfigure}[c]{0.27\textwidth}
      \caption{RadauIIA(1)}
      \begin{tikzpicture}[scale=0.5]
        \begin{axis}[xtick={8,16,32,64,128,256},
            ylabel=error,
            xlabel={$N$},
            xmode=log,
            ymode=log,
            log basis x={2},
            log basis y={10},
            ymin=1e-9,
            ymax=0.05,
            label style={font=\LARGE},
            tick label style={font=\LARGE}]
          \addplot[dashed, red, line width=2]
          table [x=N,y=H1err,col sep = comma]{heat_err_RIIA1_deg2_cfl1.csv};
          \addplot[densely dashed, blue, line width=2]
          table [x=N,y=H1err,col sep = comma]{heat_err_RIIA1_deg2_cfl4.csv};
          \addplot[loosely dashed, green, line width=2]
          table [x=N,y=H1err,col sep = comma]{heat_err_RIIA1_deg2_cfl8.csv};
          \addplot[dotted, purple, line width=2]
          table [x=N,y=L2err,col sep = comma]{heat_err_RIIA1_deg2_cfl1.csv};
          \addplot[densely dotted, orange, line width=2]
          table [x=N,y=L2err,col sep = comma]{heat_err_RIIA1_deg2_cfl4.csv};
          \addplot[loosely dotted, pink, line width=2]
          table [x=N,y=L2err,col sep = comma]{heat_err_RIIA1_deg2_cfl8.csv};
          \addplot[domain=16:64] {3/pow(x,2)} node[above, midway, yshift=-1pt, anchor=south west] {\LARGE $\mathcal{O}(h^2)$};
          \addplot[domain=16:64] {0.05/pow(x,3)} node[below, midway, xshift=-10pt] {\LARGE $\mathcal{O}(h^3)$};
        \end{axis}
      \end{tikzpicture}
    \end{subfigure}
    \begin{subfigure}[c]{0.27\textwidth}
      \caption{RadauIIA(2)}
      \begin{tikzpicture}[scale=0.5]
        \begin{axis}[xtick={8,16,32,64,128,256},
            ylabel=error,
            xlabel={$N$},
            xmode=log,
            ymode=log,
            log basis x={2},
            log basis y={10},
            ymin=1e-9,
            ymax=0.05,
            label style={font=\LARGE},
            tick label style={font=\LARGE}]
          \addplot[dashed, red, line width=2]
          table [x=N,y=H1err,col sep = comma]{heat_err_RIIA2_deg2_cfl1.csv};
          \addplot[densely dashed, blue, line width=2]
          table [x=N,y=H1err,col sep = comma]{heat_err_RIIA2_deg2_cfl4.csv};
          \addplot[loosely dashed, green, line width=2]
          table [x=N,y=H1err,col sep = comma]{heat_err_RIIA2_deg2_cfl8.csv};
          \addplot[dotted, purple, line width=2]
          table [x=N,y=L2err,col sep = comma]{heat_err_RIIA2_deg2_cfl1.csv};
          \addplot[densely dotted, orange, line width=2]
          table [x=N,y=L2err,col sep = comma]{heat_err_RIIA2_deg2_cfl4.csv};
          \addplot[loosely dotted, pink, line width=2]
          table [x=N,y=L2err,col sep = comma]{heat_err_RIIA2_deg2_cfl8.csv};
          \addplot[domain=16:64] {3/pow(x,2)} node[above, midway, yshift=-1pt, anchor=south west] {\LARGE $\mathcal{O}(h^2)$};
          \addplot[domain=16:64] {0.05/pow(x,3)} node[below, midway, xshift=-10pt] {\LARGE $\mathcal{O}(h^3)$};
        \end{axis}
      \end{tikzpicture}
    \end{subfigure}   
    \begin{subfigure}[c]{0.27\textwidth}
      \caption{RadauIIA(3)}
      \begin{tikzpicture}[scale=0.5]
        \begin{axis}[xtick={8,16,32,64,128,256},
            ylabel=error,
            xlabel={$N$},
            xmode=log,
            ymode=log,
            log basis x={2},
            log basis y={10},
            ymin=1e-9,
            ymax=0.05,
            label style={font=\LARGE},
            tick label style={font=\LARGE}]
          \addplot[dashed, red, line width=2]
          table [x=N,y=H1err,col sep = comma]{heat_err_RIIA3_deg2_cfl1.csv};
          \addplot[densely dashed, blue, line width=2]
          table [x=N,y=H1err,col sep = comma]{heat_err_RIIA3_deg2_cfl4.csv};
          \addplot[loosely dashed, green, line width=2]
          table [x=N,y=H1err,col sep = comma]{heat_err_RIIA3_deg2_cfl8.csv};
          \addplot[dotted, purple, line width=2]
          table [x=N,y=L2err,col sep = comma]{heat_err_RIIA3_deg2_cfl1.csv};
          \addplot[densely dotted, orange, line width=2]
          table [x=N,y=L2err,col sep = comma]{heat_err_RIIA3_deg2_cfl4.csv};
          \addplot[loosely dotted, pink, line width=2]
          table [x=N,y=L2err,col sep = comma]{heat_err_RIIA3_deg2_cfl8.csv};
          \addplot[domain=16:64] {3/pow(x,2)} node[above, midway, yshift=-1pt, anchor=south west] {\LARGE $\mathcal{O}(h^2)$};
          \addplot[domain=16:64] {0.05/pow(x,3)} node[below, midway, xshift=-10pt] {\LARGE $\mathcal{O}(h^3)$};
        \end{axis}
      \end{tikzpicture}
    \end{subfigure} \\
    \begin{subfigure}[c]{0.27\textwidth}
      \caption{Alexander}
      \begin{tikzpicture}[scale=0.5]
        \begin{axis}[xtick={8,16,32,64,128,256},
            ylabel=error,
            xlabel={$N$},
            xmode=log,
            ymode=log,
            log basis x={2},
            log basis y={10},
            ymin=1e-9,
            ymax=0.05,
            label style={font=\LARGE},
            tick label style={font=\LARGE}]
          \addplot[dashed, red, line width=2]
          table [x=N,y=H1err,col sep = comma]{heat_err_Alexander_deg2_cfl1.csv};
          \addplot[densely dashed, blue, line width=2]
          table [x=N,y=H1err,col sep = comma]{heat_err_Alexander_deg2_cfl4.csv};
          \addplot[loosely dashed, green, line width=2]
          table [x=N,y=H1err,col sep = comma]{heat_err_Alexander_deg2_cfl8.csv};
          \addplot[dotted, purple, line width=2]
          table [x=N,y=L2err,col sep = comma]{heat_err_Alexander_deg2_cfl1.csv};
          \addplot[densely dotted, orange, line width=2]
          table [x=N,y=L2err,col sep = comma]{heat_err_Alexander_deg2_cfl4.csv};
          \addplot[loosely dotted, pink, line width=2]
          table [x=N,y=L2err,col sep = comma]{heat_err_Alexander_deg2_cfl8.csv};
          \addplot[domain=16:64] {3/pow(x,2)} node[above, midway, yshift=-1pt, anchor=south west] {\LARGE $\mathcal{O}(h^2)$};
          \addplot[domain=16:64] {0.05/pow(x,3)} node[below, midway, xshift=-10pt] {\LARGE $\mathcal{O}(h^3)$};
        \end{axis}
      \end{tikzpicture}
    \end{subfigure}
    \begin{subfigure}[c]{0.27\textwidth}
      \caption{WSODIRK433}
      \begin{tikzpicture}[scale=0.5]
        \begin{axis}[xtick={8,16,32,64,128,256},
            ylabel=error,
            xlabel={$N$},
            xmode=log,
            ymode=log,
            log basis x={2},
            log basis y={10},
            ymin=1e-9,
            ymax=0.05,
            label style={font=\LARGE},
            tick label style={font=\LARGE}]
          \addplot[dashed, red, line width=2]
          table [x=N,y=H1err,col sep = comma]{heat_err_WSODIRK433_deg2_cfl1.csv};
          \addplot[densely dashed, blue, line width=2]
          table [x=N,y=H1err,col sep = comma]{heat_err_WSODIRK433_deg2_cfl4.csv};
          \addplot[loosely dashed, green, line width=2]
          table [x=N,y=H1err,col sep = comma]{heat_err_WSODIRK433_deg2_cfl8.csv};
          \addplot[dotted, purple, line width=2]
          table [x=N,y=L2err,col sep = comma]{heat_err_WSODIRK433_deg2_cfl1.csv};
          \addplot[densely dotted, orange, line width=2]
          table [x=N,y=L2err,col sep = comma]{heat_err_WSODIRK433_deg2_cfl4.csv};
          \addplot[loosely dotted, pink, line width=2]
          table [x=N,y=L2err,col sep = comma]{heat_err_WSODIRK433_deg2_cfl8.csv};
          \addplot[domain=16:64] {3/pow(x,2)} node[above, midway, yshift=-1pt, anchor=south west] {\LARGE $\mathcal{O}(h^2)$};
          \addplot[domain=16:64] {0.05/pow(x,3)} node[below, midway, xshift=-10pt] {\LARGE $\mathcal{O}(h^3)$};
        \end{axis}
      \end{tikzpicture}
    \end{subfigure}
    \begin{subfigure}[c]{0.27\textwidth}
      \caption{RadauIIA(4)}
      \begin{tikzpicture}[scale=0.5]
        \begin{axis}[xtick={8,16,32,64,128,256},
            ylabel=error,
            xlabel={$N$},
            xmode=log,
            ymode=log,
            log basis x={2},
            log basis y={10},
            ymin=1e-9,
            ymax=0.05,
            label style={font=\LARGE},
            tick label style={font=\LARGE},
            legend style={font=\LARGE},
            legend pos = outer north east]
          \addplot[dashed, red, line width=2]
          table [x=N,y=H1err,col sep = comma]{heat_err_RIIA4_deg2_cfl1.csv};
          \addlegendentry{$H^1$, $\Delta t=1/N$};
          \addplot[densely dashed, blue, line width=2]
          table [x=N,y=H1err,col sep = comma]{heat_err_RIIA4_deg2_cfl4.csv};
          \addlegendentry{$H^1$, $\Delta t=4/N$};
          \addplot[loosely dashed, green, line width=2]
          table [x=N,y=H1err,col sep = comma]{heat_err_RIIA4_deg2_cfl8.csv};
          \addlegendentry{$H^1$, $\Delta t=8/N$};
          \addplot[dotted, purple, line width=2]
          table [x=N,y=L2err,col sep = comma]{heat_err_RIIA4_deg2_cfl1.csv};
          \addlegendentry{$L^2$, $\Delta t=1/N$};
          \addplot[densely dotted, orange, line width=2]
          table [x=N,y=L2err,col sep = comma]{heat_err_RIIA4_deg2_cfl4.csv};
          \addlegendentry{$L^2$, $\Delta t=4/N$};
          \addplot[loosely dotted, pink, line width=2]
          table [x=N,y=L2err,col sep = comma]{heat_err_RIIA4_deg2_cfl8.csv};
          \addlegendentry{$L^2$, $\Delta t=8/N$};          
          \addplot[domain=16:64] {3/pow(x,2)} node[above, midway, yshift=-1pt, anchor=south west] {\LARGE $\mathcal{O}(h^2)$};
          \addplot[domain=16:64] {0.05/pow(x,3)} node[below, midway, xshift=-10pt] {\LARGE $\mathcal{O}(h^3)$};
        \end{axis}
      \end{tikzpicture}
    \end{subfigure}
    \caption{Accuracy vs mesh size for the heat equation on an $N \times N$ mesh of degree 2 serendipity elements.}
\label{fig:s2}
\end{center}
\end{figure}
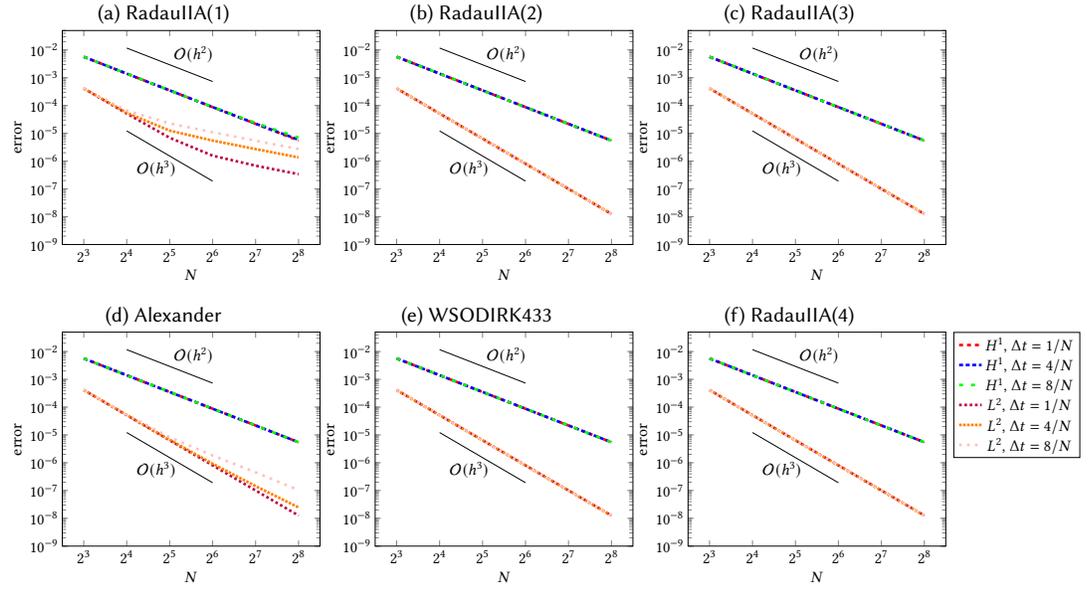

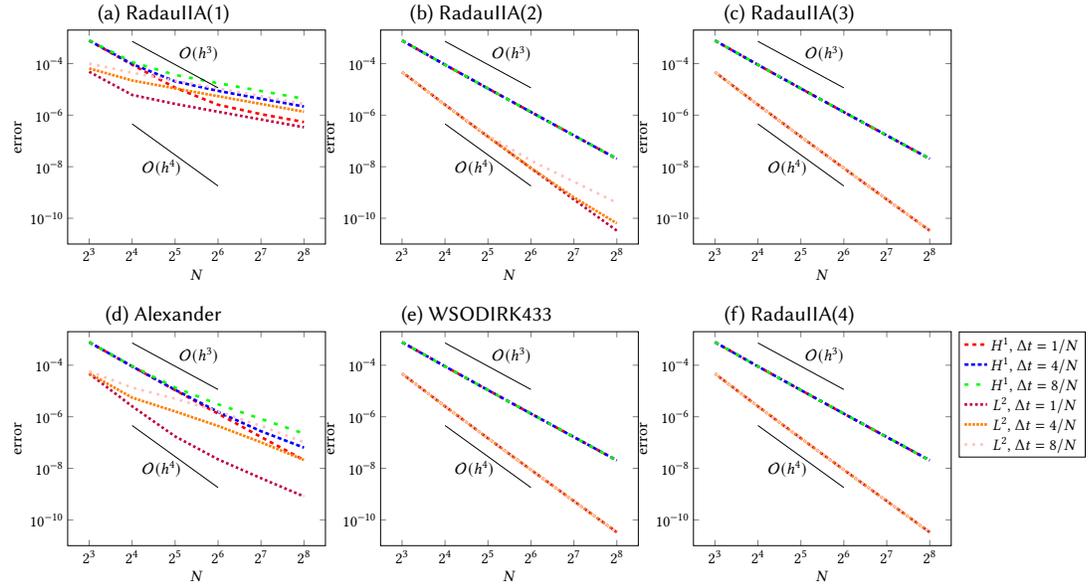
\begin{figure}
  \begin{center}
    \begin{subfigure}[c]{0.27\textwidth}
      \caption{RadauIIA(1)}
      \begin{tikzpicture}[scale=0.5]
        \begin{axis}[xtick={8,16,32,64,128,256},
            ylabel=error,
            xlabel={$N$},
            xmode=log,
            ymode=log,
            log basis x={2},
            log basis y={10},
            ymin=1e-11,
            ymax=0.002,
            label style={font=\LARGE},
            tick label style={font=\LARGE}]
          \addplot[dashed, red, line width=2]
          table [x=N,y=H1err,col sep = comma]{heat_err_RIIA1_deg3_cfl1.csv};
          \addplot[densely dashed, blue, line width=2]
          table [x=N,y=H1err,col sep = comma]{heat_err_RIIA1_deg3_cfl4.csv};
          \addplot[loosely dashed, green, line width=2]
          table [x=N,y=H1err,col sep = comma]{heat_err_RIIA1_deg3_cfl8.csv};
          \addplot[dotted, purple, line width=2]
          table [x=N,y=L2err,col sep = comma]{heat_err_RIIA1_deg3_cfl1.csv};
          \addplot[densely dotted, orange, line width=2]
          table [x=N,y=L2err,col sep = comma]{heat_err_RIIA1_deg3_cfl4.csv};
          \addplot[loosely dotted, pink, line width=2]
          table [x=N,y=L2err,col sep = comma]{heat_err_RIIA1_deg3_cfl8.csv};
          \addplot[domain=16:64] {3/pow(x,3)} node[above, midway, yshift=-1pt, anchor=south west] {\LARGE $\mathcal{O}(h^3)$};
          \addplot[domain=16:64] {0.03/pow(x,4)} node[below, midway, xshift=-10pt] {\LARGE $\mathcal{O}(h^4)$};
        \end{axis}
      \end{tikzpicture}
    \end{subfigure}
    \begin{subfigure}[c]{0.27\textwidth}
      \caption{RadauIIA(2)}
      \begin{tikzpicture}[scale=0.5]
        \begin{axis}[xtick={8,16,32,64,128,256},
            ylabel=error,
            xlabel={$N$},
            xmode=log,
            ymode=log,
            log basis x={2},
            log basis y={10},
            ymin=1e-11,
            ymax=0.002,
            label style={font=\LARGE},
            tick label style={font=\LARGE}]
          \addplot[dashed, red, line width=2]
          table [x=N,y=H1err,col sep = comma]{heat_err_RIIA2_deg3_cfl1.csv};
          \addplot[densely dashed, blue, line width=2]
          table [x=N,y=H1err,col sep = comma]{heat_err_RIIA2_deg3_cfl4.csv};
          \addplot[loosely dashed, green, line width=2]
          table [x=N,y=H1err,col sep = comma]{heat_err_RIIA2_deg3_cfl8.csv};
          \addplot[dotted, purple, line width=2]
          table [x=N,y=L2err,col sep = comma]{heat_err_RIIA2_deg3_cfl1.csv};
          \addplot[densely dotted, orange, line width=2]
          table [x=N,y=L2err,col sep = comma]{heat_err_RIIA2_deg3_cfl4.csv};
          \addplot[loosely dotted, pink, line width=2]
          table [x=N,y=L2err,col sep = comma]{heat_err_RIIA2_deg3_cfl8.csv};
          \addplot[domain=16:64] {3/pow(x,3)} node[above, midway, yshift=-1pt, anchor=south west] {\LARGE $\mathcal{O}(h^3)$};
          \addplot[domain=16:64] {0.03/pow(x,4)} node[below, midway, xshift=-10pt] {\LARGE $\mathcal{O}(h^4)$};
        \end{axis}
      \end{tikzpicture}
    \end{subfigure}    
    \begin{subfigure}[c]{0.27\textwidth}
      \caption{RadauIIA(3)}
      \begin{tikzpicture}[scale=0.5]
        \begin{axis}[xtick={8,16,32,64,128,256},
            ylabel=error,
            xlabel={$N$},
            xmode=log,
            ymode=log,
            log basis x={2},
            log basis y={10},
            ymin=1e-11,
            ymax=0.002,
            label style={font=\LARGE},
            tick label style={font=\LARGE}]
          \addplot[dashed, red, line width=2]
          table [x=N,y=H1err,col sep = comma]{heat_err_RIIA3_deg3_cfl1.csv};
          \addplot[densely dashed, blue, line width=2]
          table [x=N,y=H1err,col sep = comma]{heat_err_RIIA3_deg3_cfl4.csv};
          \addplot[loosely dashed, green, line width=2]
          table [x=N,y=H1err,col sep = comma]{heat_err_RIIA3_deg3_cfl8.csv};
          \addplot[dotted, purple, line width=2]
          table [x=N,y=L2err,col sep = comma]{heat_err_RIIA3_deg3_cfl1.csv};
          \addplot[densely dotted, orange, line width=2]
          table [x=N,y=L2err,col sep = comma]{heat_err_RIIA3_deg3_cfl4.csv};
          \addplot[loosely dotted, pink, line width=2]
          table [x=N,y=L2err,col sep = comma]{heat_err_RIIA3_deg3_cfl8.csv};
          \addplot[domain=16:64] {3/pow(x,3)} node[above, midway, yshift=-1pt, anchor=south west] {\LARGE $\mathcal{O}(h^3)$};
          \addplot[domain=16:64] {0.03/pow(x,4)} node[below, midway, xshift=-10pt] {\LARGE $\mathcal{O}(h^4)$};
        \end{axis}
      \end{tikzpicture}
    \end{subfigure} \\
    \begin{subfigure}[c]{0.27\textwidth}
      \caption{Alexander}
      \begin{tikzpicture}[scale=0.5]
        \begin{axis}[xtick={8,16,32,64,128,256},
            ylabel=error,
            xlabel={$N$},
            xmode=log,
            ymode=log,
            log basis x={2},
            log basis y={10},
            ymin=1e-11,
            ymax=0.002,
            label style={font=\LARGE},
            tick label style={font=\LARGE}]
          \addplot[dashed, red, line width=2]
          table [x=N,y=H1err,col sep = comma]{heat_err_Alexander_deg3_cfl1.csv};
          \addplot[densely dashed, blue, line width=2]
          table [x=N,y=H1err,col sep = comma]{heat_err_Alexander_deg3_cfl4.csv};
          \addplot[loosely dashed, green, line width=2]
          table [x=N,y=H1err,col sep = comma]{heat_err_Alexander_deg3_cfl8.csv};
          \addplot[dotted, purple, line width=2]
          table [x=N,y=L2err,col sep = comma]{heat_err_Alexander_deg3_cfl1.csv};
          \addplot[densely dotted, orange, line width=2]
          table [x=N,y=L2err,col sep = comma]{heat_err_Alexander_deg3_cfl4.csv};
          \addplot[loosely dotted, pink, line width=2]
          table [x=N,y=L2err,col sep = comma]{heat_err_Alexander_deg3_cfl8.csv};
          \addplot[domain=16:64] {3/pow(x,3)} node[above, midway, yshift=-1pt, anchor=south west] {\LARGE $\mathcal{O}(h^3)$};
          \addplot[domain=16:64] {0.03/pow(x,4)} node[below, midway, xshift=-10pt] {\LARGE $\mathcal{O}(h^4)$};
        \end{axis}
      \end{tikzpicture}
    \end{subfigure}
    \begin{subfigure}[c]{0.27\textwidth}
      \caption{WSODIRK433}      
      \begin{tikzpicture}[scale=0.5]
        \begin{axis}[xtick={8,16,32,64,128,256},
            ylabel=error,
            xlabel={$N$},
            xmode=log,
            ymode=log,
            log basis x={2},
            log basis y={10},
            ymin=1e-11,
            ymax=0.002,
            label style={font=\LARGE},
            tick label style={font=\LARGE}]
          \addplot[dashed, red, line width=2]
          table [x=N,y=H1err,col sep = comma]{heat_err_WSODIRK433_deg3_cfl1.csv};
          \addplot[densely dashed, blue, line width=2]
          table [x=N,y=H1err,col sep = comma]{heat_err_WSODIRK433_deg3_cfl4.csv};
          \addplot[loosely dashed, green, line width=2]
          table [x=N,y=H1err,col sep = comma]{heat_err_WSODIRK433_deg3_cfl8.csv};
          \addplot[dotted, purple, line width=2]
          table [x=N,y=L2err,col sep = comma]{heat_err_WSODIRK433_deg3_cfl1.csv};
          \addplot[densely dotted, orange, line width=2]
          table [x=N,y=L2err,col sep = comma]{heat_err_WSODIRK433_deg3_cfl4.csv};
          \addplot[loosely dotted, pink, line width=2]
          table [x=N,y=L2err,col sep = comma]{heat_err_WSODIRK433_deg3_cfl8.csv};
          \addplot[domain=16:64] {3/pow(x,3)} node[above, midway, yshift=-1pt, anchor=south west] {\LARGE $\mathcal{O}(h^3)$};
          \addplot[domain=16:64] {0.03/pow(x,4)} node[below, midway, xshift=-10pt] {\LARGE $\mathcal{O}(h^4)$};
        \end{axis}
      \end{tikzpicture}
    \end{subfigure}     
    \begin{subfigure}[c]{0.27\textwidth}
      \caption{RadauIIA(4)}
      \begin{tikzpicture}[scale=0.5]
        \begin{axis}[xtick={8,16,32,64,128,256},
            ylabel=error,
            xlabel={$N$},
            xmode=log,
            ymode=log,
            log basis x={2},
            log basis y={10},
            ymin=1e-11,
            ymax=0.002,
            label style={font=\LARGE},
            tick label style={font=\LARGE},
            legend style={font=\LARGE},
            legend pos = outer north east]
          \addplot[dashed, red, line width=2]
          table [x=N,y=H1err,col sep = comma]{heat_err_RIIA4_deg3_cfl1.csv};
          \addlegendentry{$H^1$, $\Delta t=1/N$};
          \addplot[densely dashed, blue, line width=2]
          table [x=N,y=H1err,col sep = comma]{heat_err_RIIA4_deg3_cfl4.csv};
          \addlegendentry{$H^1$, $\Delta t=4/N$};
          \addplot[loosely dashed, green, line width=2]
          table [x=N,y=H1err,col sep = comma]{heat_err_RIIA4_deg3_cfl8.csv};
          \addlegendentry{$H^1$, $\Delta t=8/N$};
          \addplot[dotted, purple, line width=2]
          table [x=N,y=L2err,col sep = comma]{heat_err_RIIA4_deg3_cfl1.csv};
          \addlegendentry{$L^2$, $\Delta t=1/N$};
          \addplot[densely dotted, orange, line width=2]
          table [x=N,y=L2err,col sep = comma]{heat_err_RIIA4_deg3_cfl4.csv};
          \addlegendentry{$L^2$, $\Delta t=4/N$};
          \addplot[loosely dotted, pink, line width=2]
          table [x=N,y=L2err,col sep = comma]{heat_err_RIIA4_deg3_cfl8.csv};
          \addlegendentry{$L^2$, $\Delta t=8/N$};
           \addplot[domain=16:64] {3/pow(x,3)} node[above, midway, yshift=-1pt, anchor=south west] {\LARGE $\mathcal{O}(h^3)$};
          \addplot[domain=16:64] {0.03/pow(x,4)} node[below, midway, xshift=-10pt] {\LARGE $\mathcal{O}(h^4)$};
        \end{axis}
      \end{tikzpicture}
    \end{subfigure}
    \caption{Accuracy vs mesh size for the heat equation on an $N \times N$ mesh of degree 3 serendipity elements.}
\label{fig:s3}
\end{center}
\end{figure}

Figures~\ref{fig:s2} and~\ref{fig:s3} plot the $H^1$ and $L^2$ errors obtained in integrating the heat equation until $T=1$ using various schemes.
On an $N \times N$ mesh, we take the time step as $\Delta t = M / N$, where $M = 1, 4, 8$.
This allows us to explore the onset of order reduction in our time stepping methods.
In Figure~\ref{fig:s2}, where we use quadratic serendipity elements,
we see suboptimal accuracy in $L^2$ for the 1-stage RadauIIA method (backward Euler) with all values of $M$.
We also see suboptimal $L^2$ accuracy with the 3-stage DIRK of Alexander using $M=8$.
For backward Euler, we only expect first-order method.  The Alexander DIRK is formally higher-order, but the very large time step creates order reduction due to increased stiffness.
Otherwise, for all higher-order RadauIIA methods and the 4-stage DIRK, we obtain optimal accuracy even with large time steps, indicating that spatial rather than temporal error dominates.
Hence, one should choose between RadauIIA($s$) with $s \geq 2$ and the four-stage DIRK as indicated by performance or some other metric.
The three-stage Alexander DIRK is possibly usable if the time step is not too large.

Figure~\ref{fig:s3} repeates this experiment with cubic serendipity elements as the spatial discretization.
Since errors should decrease more quickly than for quadratic elements, this experiment should be more sensitive to a loss in accuracy from the time-stepper.
Unsurprisingly, RadauIIA(1) gives poor results. The loss of accuracy for the three-stage DIRK is more pronounced now.
We also see some limits in accuracy using the RadauIIA(2) method, which is formally third order (with stage order two).
No loss of accuracy appears in the $H^1$ error, but the $L^2$ error suffers with large time steps.
The higher-order RadauIIA methods and the four-stage DIRK here all give comparable accuracy, with temporal integration error apparently dominated by spatial discretization.

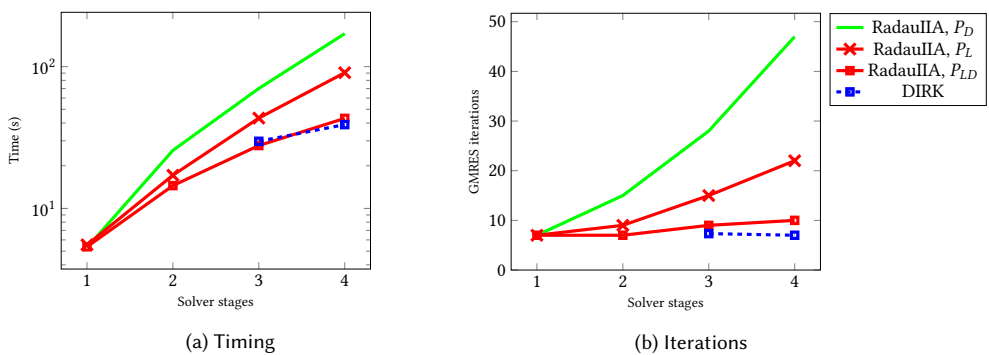
\begin{figure}
  \begin{subfigure}[c]{0.40\textwidth}
    \begin{tikzpicture}[scale=0.6]
      \begin{axis}[ymode=log,
          xtick={1,2,3,4},
          ylabel={Time (s)},
          xlabel={Solver stages},
          tick label style={font=\LARGE}]
        \addplot[green, line width = 2]
        table [x=ns,y=time,col sep=comma] {heat_perf_bjac_firk.csv};
        \addplot[red, mark=x, mark size = 5, line width = 2]
        table [x=ns,y=time,col sep=comma] {heat_perf_bgs_firk.csv};
        \addplot[red, mark=square, line width = 2]
        table [x=ns,y=time,col sep=comma] {heat_perf_rana_amg_firk.csv};
        \addplot[blue, dashed, mark=square, mark options = {solid}, line width = 2]
        table [x=ns,y=time,col sep=comma] {heat_perf_amg_dirk.csv};        
      \end{axis}
    \end{tikzpicture}
    \caption{Timing}
    \label{fig:heatperftime}
  \end{subfigure}
  \begin{subfigure}[c]{0.40\textwidth}
    \begin{tikzpicture}[scale=0.6]
      \begin{axis}[xtick={1,2,3,4}, ymin=0,
        legend style={font=\LARGE},
        legend pos = outer north east,
        ylabel={GMRES iterations},
        xlabel={Solver stages},
            tick label style={font=\LARGE}]
        \addplot[green, line width = 2]
        table [x=ns,y=Its,col sep=comma] {heat_perf_bjac_firk.csv};
        \addlegendentry{RadauIIA, $P_{D}$};
        \addplot[red, mark=x, mark size = 5, line width = 2]
        table [x=ns,y=Its,col sep=comma] {heat_perf_bgs_firk.csv};
        \addlegendentry{RadauIIA, $P_{L}$};
        \addplot[red, mark=square, line width = 2]
        table [x=ns,y=Its,col sep=comma] {heat_perf_rana_amg_firk.csv};
        \addlegendentry{RadauIIA, $P_{LD}$};
        \addplot[blue, dashed, mark=square, mark options = {solid}, line width = 2]
        table [x=ns,y=Its,col sep=comma] {heat_perf_amg_dirk.csv};
        \addlegendentry{DIRK};        
      \end{axis}[dashed, line width = 2]
    \end{tikzpicture}
    \caption{Iterations}
    \label{fig:heatperfits}    
  \end{subfigure}
  \caption{Solver performance as a function of stage count for the heat equation.  We integrated the heat equation using cubic serendipity elements on a $64 \times 64$ mesh using various solver configurations for RadauIIA and DIRK methods.}
  \label{fig:heatperf}
\end{figure}

Now, we consider performance of these methods with various solver strategies.
For the heat equation, \lstinline{hypre}~\cite{falgout2002hypre} provides a highly capable preconditioner for single-stage methods, and so we focus on stage-segregated preconditioners such as~\eqref{eq:pcbjac}, ~\eqref{eq:pcbgs}, and~\eqref{eq:pcranald} for the fully implicit methods and compare the performance to that of diagonally implicit methods.
We compare the total run-time and Krylov iterations per time step a fixed $256 \times 256$ mesh, integrating over a time interval $[0, 2]$ with the time steppers above and various solution techniques.
In each case, we solve the linear system at each time step with FGMRES with residual relative tolerance of $10^{-8}$.  The restart threshhold is set at 50, which is well above the number of iterations used in any of our computations.
For the multi-stage RadauIIA methods, we consider the block Jacobi~\eqref{eq:pcbjac}, block Gauss-Seidel~\eqref{eq:pcbgs}, and Rana-LD~\eqref{eq:pcranald} preconditioners, with the inverse of diagonal blocks approximated by a single application of \lstinline{hypre} with default parameters.  For the single-stage RadauIIA and each stage of the diagonally implicit methods, we simply use FGMRES preconditioned with a single sweep of \lstinline{hypre} to solve each stage.

Figure~\ref{fig:heatperf} shows the solver performance as a function of stage count.  On the right, we report the average number of Krylov iterations per linear solve.  For fully implicit methods, this is equivalent to the number of Krylov iterations per time step but, for diagonally implicit methods, we are averaging over the stages as well (reporting total number of Krylov iterations divided by the product of the number of stages and the number of time steps).
Here, we see that the diagonally implicit method requires slightly fewer Krylov iterations than the Rana-type preconditioner~$P_{LD}$ from~\eqref{eq:pcranald}, although the iteration count is fairly flat over the number of stages.  On the other hand, the block Jacobi preconditioner~$P_{D}$ from~\eqref{eq:pcbjac} and block Gauss-Seidel preconditioner~$P_{L}$ from~\eqref{eq:pcbgs} give iteration counts growing with the number of stages.
Figure~\ref{fig:heatperftime} shows the total run-time for the time-stepping loop.  This includes assembling matrices and setting up PETSc solvers as well as solving linear systems and updating the solution at each time step.  Here, we see a performance gap between the block Jacobi and Gauss-Seidel preconditioners on one hand and the diagonally implicit and Rana-type preconditioners on the other.
However, the run-time is very comparable between the two latter methods, even though DIRK methods require slightly fewer iterations.

\subsection{Navier-Stokes}
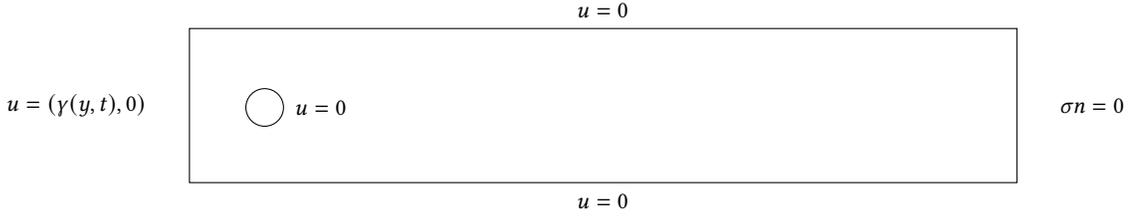
\begin{figure}
  \begin{center}
  \begin{tikzpicture}[scale=5]
    \draw (0,0) rectangle (2.2, 0.41);
    \draw (.2, .2) circle (0.05);
    \node[] at (1.1, 0.46) {$u=0$};
    \node[] at (1.1, -0.05) {$u=0$};
    \node[] at (.35, .2) {$u=0$};
    \node[] at (2.4, 0.205) {$\sigma n = 0$};
    \node[] at (-0.3, 0.205) {$u=\left(\gamma(y, t), 0\right)$};
  \end{tikzpicture}
  \end{center}
  \caption{Computational domain for Navier-Stokes flow past a cylinder}
  \Description{Picture showing computational domain for flow past a cylinder.}
  \label{fig:domain}
\end{figure}

To illustrate the efficacy of our techniques for the Navier-Stokes equations, 
we consider a common benchmark of computing the drag and lift on a flow past a cylinder.  
The domain, shown in Figure~\ref{fig:domain}, consists of the rectangle $[0, 2.2] \times [0, 0.41]$ with the circle of radius 0.05 centered at (0.2, 0.2) omitted from the domain.  
On the left edge, we impose an inflow condition, setting the horizontal velocity component to be
\[ 
\gamma(y, t) = 6 \sin \left(\tfrac{\pi t}{8} \right) \left( \tfrac{y(0.41-y)}{0.41^2} \right)
\]
and the vertical velocity component to be zero.  
Along the top and bottom edges and boundary of the circle, we impose no-slip conditions, and we have natural outflow conditions on the right edge.

Using Firedrake's interface to netgen~\cite{betteridge2024ngspetsc}, we built a mesh hierarchy with three levels, the finest mesh comprising 12,800 cells and 6,616 vertices.  With each refinement, vertices introduced at the midpoints of edges along the curved boundary of the circle are projected onto the circle, to improve the geometric resolution of the cylinder.

We integrate the system until $T=8$, measuring the drag and lift on the circle boundary at each time step, as in~\cite{john2004reference,farrell2021irksome,abu2022monolithic}.   
Our goal is to study the efficiency and accuracy of time-stepping for this problem.  
Hence, we fix the spatial discretization and vary the Runge-Kutta method and time step size.
Since the standard Taylor-Hood methods lose inf-sup stability when the time step is decreased on a fixed mesh, we use the Mardal-Tai-Winther finite element for velocities paired with piecewise constant pressures~\cite{mardal2002robust}.
This pair is stable and conforming for the mixed Poisson problem on $\hdiv \times L^2$ and a suitable nonconforming pair on $(H^1)^2 \times L^2$.
We refer the reader to~\cite{farrell2022transformations} for details on the inclusion of the MTW element in Firedrake.

For each Runge-Kutta method, we compute each time-step using a Newton-Krylov method.
Convergence for each time step is measured by requiring a reduction in either the absolute Euclidean norm of the nonlinear residual for the stage-coupled system below $10^{-12}$ or a relative reduction in the nonlinear residual norm by the same factor, with the FGMRES convergence tolerance for each linearization chosen adaptively according to the Eisenstat-Walker technique~\cite{eisenstat1996choosing}.
In these cases, we found that the solver required somewhat fewer iterations and hence gave better run-time with the \lstinline{IA} splitting and stage value formulation, so we report using this configuration.

\begin{figure}[h]
  \centering
  \begin{tikzpicture}[scale=0.8]
        \draw[thin] (5,5) -- (9,5) -- (9, 9) -- (5, 9) -- (5,5);
        \draw[thin] (7,5) -- (5,7);
        \draw[thin] (9,5) -- (7,7);
        \draw[thin] (7,7) -- (5,9);
        \draw[thin] (9,7) -- (7,9);
        \draw[thin] (9,7) -- (5,7);
        \draw[thin] (7,5) -- (7,9);
        \draw[very thick,->] (7.5, 7) -- (7.5, 7-.4);
        \draw[very thick,->] (8.5, 7) -- (8.5, 7-.4);
        \draw[very thick,->] (7.5, 8.5) -- (7.5+.283, 8.5+.283);
        \draw[very thick,->] (8.5, 7.5) -- (8.5+.283, 7.5+.283);
        \draw[very thick, ->] (7, 7.5) -- (7-.4, 7.5);
        \draw[very thick, ->] (7, 8.5) -- (7-.4, 8.5);
        \draw[very thick, ->] (7.8, 7) -- (8.2, 7);
        \draw[very thick, ->] (7, 7.8) -- (7, 8.2);
        \draw[very thick, ->] (7.5 + .359, 8.5 - .359) -- (8.5 - .359, 7.5 + .359);
        \node[circle,fill=black,inner sep=0pt,minimum size=7pt] at (7.65,7.7) {};
\end{tikzpicture}
\caption{Typical Vanka patch for MTW with arrows showing velocity degrees of freedom and pressure as a black circle.}
\label{fig:vanka}
\end{figure}
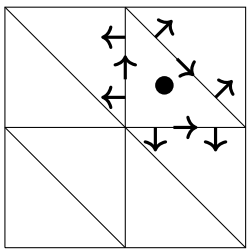

Typically, only a small handful of Newton iterations were used, and we never approached the PETSc default maximum iteration count of 50.
We set a maximum Krylov iteration count of 200, but never approached this, either.
We considered two preconditioners for the Jacobian in fully implicit methods.
First, we considered a monolithic multigrid scheme in which
FGMRES is preconditioned with geometric multigrid V(3,3) cycles using elementwise Vanka relaxation on each level but the coarsest, GMRES to accelerate relaxation on each level, and a sparse direct method on the coarsest mesh.
The typical Vanka patch is shown in Figure~\ref{fig:vanka}, where we take the single pressure degree of freedom within an element together with all velocity degrees of freedom associated with that element.
This is implemented through Firedrake's \lstinline{ASMVankaPC}.
In the multi-stage case, we use monolithic multigrid that includes all of the velocity and pressure degrees of freedom from all RK stages in a single patch, as done in~\cite{abu2022monolithic}.
Additionally, we considered a stage-segregated Rana-type preconditioner~\eqref{eq:pcranald}, in which we approximate the inverse of the diagonal blocks with the single-stage version of the multigrid scheme described above.  We note that using ``black-box'' numerical solvers, such as AMG, is not possible for this discretization (nor for that of Cahn-Hilliard, discussed below), so we can offer no comparison to previously available solvers for this problem.

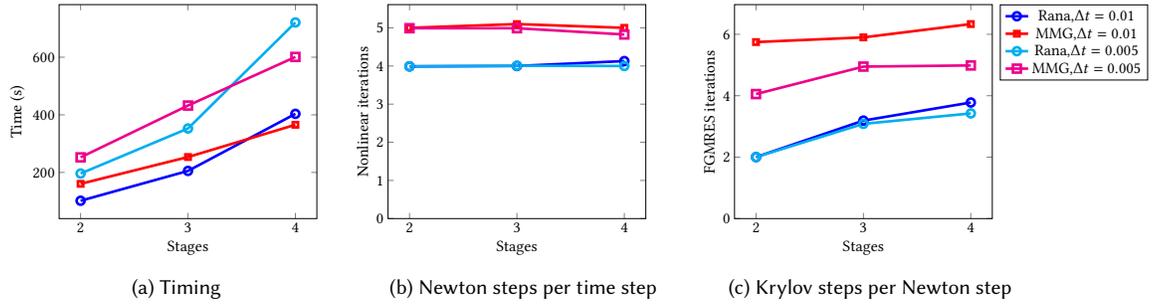
\begin{figure}
  \begin{center}
    \begin{subfigure}[c]{0.3\textwidth}
    \begin{tikzpicture}[scale=0.5]
      \begin{axis}[xlabel=Stages, ylabel={Time (s)},
          xtick={2,3,4},
          label style={font=\LARGE},
          tick label style={font=\LARGE}]
        \addplot[blue, line width = 2, mark=o, mark size=3]
        table [x=Stages,y=time, col sep=comma]{RIIA.rana.dt0.01.csv};
        \addplot[red, mark=square, mark options = {solid}, line width = 2]
        table [x=Stages,y=time, col sep=comma]{RIIA.vanka.dt0.01.csv};
        \addplot[cyan, mark=o, mark size = 3, mark options = {solid}, line width = 2]
        table [x=Stages,y=time, col sep=comma]{RIIA.rana.dt0.005.csv};
        \addplot[magenta, mark=square, mark size = 3, mark options = {solid}, line width = 2]
        table [x=Stages,y=time, col sep=comma]{RIIA.vanka.dt0.005.csv};
      \end{axis}
    \end{tikzpicture}
    \caption{Timing}
    \label{fig:ranavmmgtime}
    \end{subfigure}
    \begin{subfigure}[c]{0.3\textwidth}
    \begin{tikzpicture}[scale=0.5]
      \begin{axis}[xlabel=Stages, ylabel={Nonlinear iterations},
          xtick={2,3,4},
          ymin=0,
          label style={font=\LARGE},
          tick label style={font=\LARGE}]
        \addplot[blue, line width = 2, mark=o, mark size=3]
        table [x=Stages,y=newt, col sep=comma]{RIIA.rana.dt0.01.csv};
        \addplot[red, mark=square, mark options = {solid}, line width = 2]
        table [x=Stages,y=newt, col sep=comma]{RIIA.vanka.dt0.01.csv};
        \addplot[cyan, mark=o, mark size = 3, mark options = {solid}, line width = 2]
        table [x=Stages,y=newt, col sep=comma]{RIIA.rana.dt0.005.csv};
        \addplot[magenta, mark=square, mark size = 3, mark options = {solid}, line width = 2]
        table [x=Stages,y=newt, col sep=comma]{RIIA.vanka.dt0.005.csv};
      \end{axis}
    \end{tikzpicture}
    \caption{Newton steps per time step}
    \label{fig:ranavmmgnewt}
    \end{subfigure}
    \begin{subfigure}[c]{0.3\textwidth}
    \begin{tikzpicture}[scale=0.5]
      \begin{axis}[xlabel=Stages, ylabel={FGMRES iterations},
          xtick={2,3,4},
          ymin=0,
          label style={font=\LARGE},
          tick label style={font=\LARGE},
          legend style={font=\LARGE},
            legend pos = outer north east]
        \addplot[blue, line width = 2, mark=o, mark size=3]
        table [x=Stages,y=ksp, col sep=comma]{RIIA.rana.dt0.01.csv};
        \addlegendentry{Rana,$\Delta t=0.01$};
        \addplot[red, mark=square, mark options = {solid}, line width = 2]
        table [x=Stages,y=ksp, col sep=comma]{RIIA.vanka.dt0.01.csv};
        \addlegendentry{MMG,$\Delta t=0.01$};        
        \addplot[cyan, mark=o, mark size = 3, mark options = {solid}, line width = 2]
        table [x=Stages,y=ksp, col sep=comma]{RIIA.rana.dt0.005.csv};
        \addlegendentry{Rana,$\Delta t=0.005$};
        \addplot[magenta, mark=square, mark size = 3, mark options = {solid}, line width = 2]
        table [x=Stages,y=ksp, col sep=comma]{RIIA.vanka.dt0.005.csv};
        \addlegendentry{MMG,$\Delta t=0.005$};        
      \end{axis}
    \end{tikzpicture}
    \caption{Krylov steps per Newton step}
    \label{fig:ranavmmggmres}
    \end{subfigure}    
  \end{center}
  \caption{Solver performance comparing Rana and monolithic multigrid (MMG) preconditioning strategies for integrating the Navier--Stokes equations on the time interval [0, 1].  The Rana scheme requires slightly fewer linear and nonlinear iterations than the monolithic multigrid scheme, but gives somewhat higher run-time with a large stage count.}
  \label{fig:ranavmmg}
\end{figure}

Before performing a full suite of benchmark calculations over the time interval [0,8], we first consider the relative performance of monolithic and Rana-type preconditioners over a smaller but still significant time interval [0, 1], which we present in Figure~\ref{fig:ranavmmg}.  We consider two different time steps, of size $0.01$ and $0.005$.
Here, we see that the Rana preconditioner is leading to somewhat fewer Newton iterations per time step (recall that Eisenstat-Walker uses an initially loose Krylov tolarance) and also fewer Krylov steps per Newton step.  This leads to a lower overall run-time for the 2- and 3-stage RadauIIA method.  However, we see that the monolithic multigrid scheme actually gives a somehwat lower run-time for RadauIIA(4).  We expect this is because as the monolithic scheme has many fewer traversals of the multigrid hierarchy and a greater arithmetic intensity in applying the relaxation.
Because it is somewhat better at high stage counts, we proceed with our full benchmark experiments using a monolithic multigrid scheme for our fully implicit methods.

Our relatively coarse spatial discretization captures the benchmark drag and lift values to a few decimal places, but our main focus here is the accuracy and efficiency of time integration.
So, we will take our reference lift/drag values to be those computed with this spatial discretization using exact time integration, which we approximate in practice with a 5-stage RadauII approximation and the finest time step.
Our reference values for the maximum drag value is approximately 2.89 at time $3.93625$ and the maximum lift value is about 0.47 at time 5.68875.

In our experiments, we integrated the system using RadauIIA methods with one through four stages and a suite of DIRK methods with various numbers of stages.
In particular, we used the method of Alexander together with
WSODIRK433 and WSODIRK744 methods.
We varied varied the time step from $\Delta t = 0.16$ down to $\Delta t = 0.0025$.  
In each case, we measured the total wall-clock time for setting up the time stepper and integrating in time and accumulated the total number of Newton and GMRES iterations performed and tracked the lift and drag coefficients at each time step.
This allows us to report the error in maximum lift and drag computed, as well as the error at $T=8$.

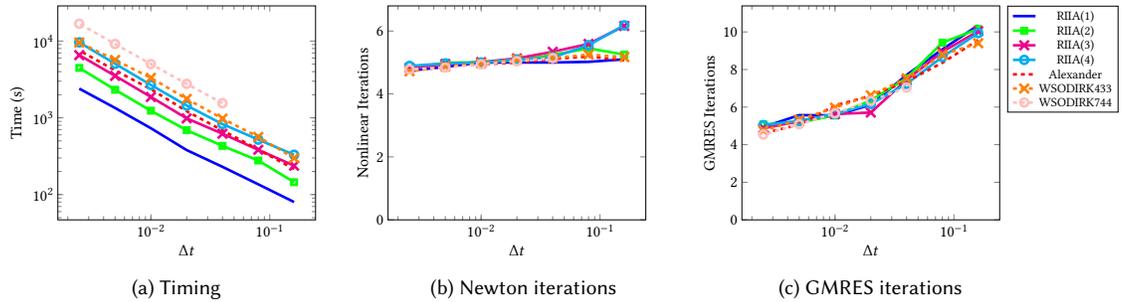
\begin{figure}
  \begin{center}
    \begin{subfigure}[c]{0.3\textwidth}
    \begin{tikzpicture}[scale=0.5]
      \begin{axis}[xlabel=$\Delta t$, ylabel={Time (s)},
          xmode=log, ymode=log,
          label style={font=\LARGE},
            tick label style={font=\LARGE}]
        \addplot[blue, line width = 2]
        table [x=dt,y=time, col sep=comma]{NSE_RadauIIA1.csv};
        \addplot[green, mark=square, mark options = {solid}, line width = 2]
        table [x=dt,y=time, col sep=comma]{NSE_RadauIIA2.csv};
        \addplot[magenta, mark=x, mark size = 5, mark options = {solid}, line width = 2]
        table [x=dt,y=time, col sep=comma]{NSE_RadauIIA3.csv};
        \addplot[cyan, mark=o, mark size = 3, mark options = {solid}, line width = 2]
        table [x=dt,y=time, col sep=comma]{NSE_RadauIIA4.csv};
        \addplot[red, dashed, line width = 2]
        table [x=dt,y=time, col sep=comma]{NSE_Alexander.csv};
        \addplot[orange, mark=x, mark size = 5, mark options = {solid}, dashed, line width = 2]
        table [x=dt,y=time, col sep=comma]{NSE_WSODIRK433.csv};
        \addplot[pink, mark=o, mark size = 3, mark options = {solid}, dashed, line width = 2]
        table [x=dt,y=time, col sep=comma]{NSE_WSODIRK744.csv};
      \end{axis}
    \end{tikzpicture}
    \caption{Timing}
    \label{fig:nsetime}
    \end{subfigure}
    \begin{subfigure}[c]{0.3\textwidth}
      \begin{tikzpicture}[scale=0.5]
        \begin{axis}[xlabel=$\Delta t$, ylabel={Nonlinear Iterations}, ymin=0, xmode=log,
            label style={font=\LARGE},
            tick label style={font=\LARGE}]
          \addplot[blue, line width = 2]
          table [x=dt,y=newt, col sep=comma]{NSE_RadauIIA1.csv};
          \addplot[green, mark=square, mark options = {solid}, line width = 2]
          table [x=dt,y=newt, col sep=comma]{NSE_RadauIIA2.csv};
          \addplot[magenta, mark=x, mark size = 5, mark options = {solid}, line width = 2]
          table [x=dt,y=newt, col sep=comma]{NSE_RadauIIA3.csv};
          \addplot[cyan, mark=o, mark size = 3, mark options = {solid}, line width = 2]
          table [x=dt,y=newt, col sep=comma]{NSE_RadauIIA4.csv};
          \addplot[red, dashed, line width = 2]
          table [x=dt,y=newt, col sep=comma]{NSE_Alexander.csv};
          \addplot[orange, mark=x, mark size = 5, mark options = {solid}, dashed, line width = 2]
          table [x=dt,y=newt, col sep=comma]{NSE_WSODIRK433.csv};
          \addplot[pink, mark=o, mark size = 3, mark options = {solid}, dashed, line width = 2]
          table [x=dt,y=newt, col sep=comma]{NSE_WSODIRK744.csv};
        \end{axis}
      \end{tikzpicture}
      \caption{Newton iterations}
      \label{fig:nsenewt}
    \end{subfigure}
  \begin{subfigure}[c]{0.3\textwidth}
    \begin{tikzpicture}[scale=0.5]
      \begin{axis}[xlabel=$\Delta t$, ylabel={GMRES Iterations}, ymin=0,xmode=log, legend style={nodes={scale=0.8},font=\LARGE},
        legend pos = outer north east,
            label style={font=\LARGE},
            tick label style={font=\LARGE}]
        \addplot[blue, line width = 2]
          table [x=dt,y=ksp, col sep=comma]{NSE_RadauIIA1.csv};
          \addlegendentry{RIIA(1)}
          \addplot[green, mark=square, mark options = {solid}, line width = 2]
          table [x=dt,y=ksp, col sep=comma]{NSE_RadauIIA2.csv};
          \addlegendentry{RIIA(2)}
          \addplot[magenta, mark=x, mark size = 5, mark options = {solid}, line width = 2]
          table [x=dt,y=ksp, col sep=comma]{NSE_RadauIIA3.csv};
          \addlegendentry{RIIA(3)}
          \addplot[cyan, mark=o, mark size = 3, mark options = {solid}, line width = 2]
          table [x=dt,y=ksp, col sep=comma]{NSE_RadauIIA4.csv};
          \addlegendentry{RIIA(4)}
          \addplot[red, dashed, line width = 2]
          table [x=dt,y=ksp, col sep=comma]{NSE_Alexander.csv};
          \addlegendentry{Alexander}
          \addplot[orange, mark=x, mark size = 5, mark options = {solid}, dashed, line width = 2]
          table [x=dt,y=ksp, col sep=comma]{NSE_WSODIRK433.csv};
          \addlegendentry{WSODIRK433}
          \addplot[pink, mark=o, mark size = 3, mark options = {solid}, dashed, line width = 2]
          table [x=dt,y=ksp, col sep=comma]{NSE_WSODIRK744.csv};
          \addlegendentry{WSODIRK744}
        \end{axis}
    \end{tikzpicture}
    \caption{GMRES iterations}
    \label{fig:nsegmres}
  \end{subfigure}
  \end{center}
  \caption{Solver performance for integrating Navier--Stokes equations on $[0,8]$.  We report total wall-clock time for the time stepping loop as well as average nonlinear iterations per time step and average linear iterations per nonlinear iteration.  For DIRKs, our Newton iterations are averaged over each nonlinear stage of each time step.}
  \label{fig:perfnse}
\end{figure}

\begin{figure}
  \begin{center}
    \begin{subfigure}[c]{0.45\textwidth}
      \begin{tikzpicture}[scale=0.8]
        \begin{axis}[xlabel=$\Delta t$, ylabel={$\epsilon_D^{\max}$}, xmode=log,ymode=log, legend pos=south east, legend style={nodes={scale=0.8}}]
          \addplot[blue, line width = 2]
          table [x=dt,y=edmax, col sep=comma]{NSE_RadauIIA1.csv};
          \addlegendentry{RIIA(1)}
          \addplot[green, mark=square, mark options = {solid}, line width = 2]
          table [x=dt,y=edmax, col sep=comma]{NSE_RadauIIA2.csv};
          \addlegendentry{RIIA(2)}
          \addplot[magenta, mark=x, mark size = 5, mark options = {solid}, line width = 2]
          table [x=dt,y=edmax, col sep=comma]{NSE_RadauIIA3.csv};
          \addlegendentry{RIIA(3)}
          \addplot[cyan, mark=o, mark size = 3, mark options = {solid}, line width = 2]
          table [x=dt,y=edmax, col sep=comma]{NSE_RadauIIA4.csv};
          \addlegendentry{RIIA(4)}
          \addplot[red, dashed, line width = 2]
          table [x=dt,y=edmax, col sep=comma]{NSE_Alexander.csv};
          \addlegendentry{Alexander}
          \addplot[orange, mark=x, mark size = 5, mark options = {solid}, dashed, line width = 2]
          table [x=dt,y=edmax, col sep=comma]{NSE_WSODIRK433.csv};
          \addlegendentry{WSODIRK433}
          \addplot[pink, mark=o, mark size = 3, mark options = {solid}, dashed, line width = 2]
          table [x=dt,y=edmax, col sep=comma]{NSE_WSODIRK744.csv};
          \addlegendentry{WSODIRK744}
          \addplot[domain=0.002:0.02] {20/pow(x,-3)} node[below, midway, yshift=2pt, anchor=east]{$\mathcal{O}(\Delta t^3)$};
        \end{axis}
      \end{tikzpicture}
      \caption{Error in maximum drag}
      \label{fig:nseedmax}
    \end{subfigure}
    \begin{subfigure}[c]{0.45\textwidth}
      \begin{tikzpicture}[scale=0.8]
        \begin{axis}[xlabel=$\Delta t$, ylabel={$\epsilon_L^{\max}$},xmode=log, ymode=log]
        \addplot[blue, line width = 2]
          table [x=dt,y=elmax, col sep=comma]{NSE_RadauIIA1.csv};
          \addplot[green, mark=square, mark options = {solid}, line width = 2]
          table [x=dt,y=elmax, col sep=comma]{NSE_RadauIIA2.csv};
          \addplot[magenta, mark=x, mark size = 5, mark options = {solid}, line width = 2]
          table [x=dt,y=elmax, col sep=comma]{NSE_RadauIIA3.csv};
          \addplot[cyan, mark=o, mark size = 3, mark options = {solid}, line width = 2]
          table [x=dt,y=elmax, col sep=comma]{NSE_RadauIIA4.csv};
          \addplot[red, dashed, line width = 2]
          table [x=dt,y=elmax, col sep=comma]{NSE_Alexander.csv};
          \addplot[orange, mark=x, mark size = 5, mark options = {solid}, dashed, line width = 2]
          table [x=dt,y=elmax, col sep=comma]{NSE_WSODIRK433.csv};
          \addplot[pink, mark=o, mark size = 3, mark options = {solid}, dashed, line width = 2]
          table [x=dt,y=elmax, col sep=comma]{NSE_WSODIRK744.csv};
          \addplot[domain=0.002:0.02] {8000/pow(x,-3)} node[below, midway, yshift=2pt, anchor=east]{$\mathcal{O}(\Delta t^3)$};
        \end{axis}
      \end{tikzpicture}
      \caption{Error in maximum lift}
      \label{fig:nseelmax}      
    \end{subfigure} \\[12pt]
    \begin{subfigure}[c]{0.45\textwidth}
      \begin{tikzpicture}[scale=0.8]
        \begin{axis}[xlabel=$\Delta t$, ylabel={$\epsilon_D(T)$}, xmode=log,ymode=log]
          \addplot[blue, line width = 2]
          table [x=dt,y=edt, col sep=comma]{NSE_RadauIIA1.csv};
          \addplot[green, mark=square, mark options = {solid}, line width = 2]
          table [x=dt,y=edt, col sep=comma]{NSE_RadauIIA2.csv};
          \addplot[magenta, mark=x, mark size = 5, mark options = {solid}, line width = 2]
          table [x=dt,y=edt, col sep=comma]{NSE_RadauIIA3.csv};
          \addplot[cyan, mark=o, mark size = 3, mark options = {solid}, line width = 2]
          table [x=dt,y=edt, col sep=comma]{NSE_RadauIIA4.csv};
          \addplot[red, dashed, line width = 2]
          table [x=dt,y=edt, col sep=comma]{NSE_Alexander.csv};
          \addplot[orange, mark=x, mark size = 5, mark options = {solid}, dashed, line width = 2]
          table [x=dt,y=edt, col sep=comma]{NSE_WSODIRK433.csv};
          \addplot[pink, mark=o, mark size = 3, mark options = {solid}, dashed, line width = 2]
          table [x=dt,y=edt, col sep=comma]{NSE_WSODIRK744.csv};
          \addplot[domain=0.002:0.02] {10/pow(x,-5)} node[below, midway, yshift=-1pt, anchor=west]{$\mathcal{O}(\Delta t^5)$};
          \addplot[domain=0.002:0.02] {500/pow(x,-3)} node[below, midway, yshift=2pt, anchor=east]{$\mathcal{O}(\Delta t^3)$};          
        \end{axis}
      \end{tikzpicture}
      \caption{Error in final drag}
      \label{fig:nseedt}     
    \end{subfigure}
    \begin{subfigure}[c]{0.45\textwidth}
      \begin{tikzpicture}[scale=0.8]
        \begin{axis}[xlabel=$\Delta t$, ylabel={$\epsilon_L(T)$},xmode=log, ymode=log, legend pos=south east]
        \addplot[blue, line width = 2]
          table [x=dt,y=elt, col sep=comma]{NSE_RadauIIA1.csv};
          \addplot[green, mark=square, mark options = {solid}, line width = 2]
          table [x=dt,y=elt, col sep=comma]{NSE_RadauIIA2.csv};
          \addplot[magenta, mark=x, mark size = 5, mark options = {solid}, line width = 2]
          table [x=dt,y=elt, col sep=comma]{NSE_RadauIIA3.csv};
          \addplot[cyan, mark=o, mark size = 3, mark options = {solid}, line width = 2]
          table [x=dt,y=elt, col sep=comma]{NSE_RadauIIA4.csv};
          \addplot[red, dashed, line width = 2]
          table [x=dt,y=elt, col sep=comma]{NSE_Alexander.csv};
          \addplot[orange, mark=x, mark size = 5, mark options = {solid}, dashed, line width = 2]
          table [x=dt,y=elt, col sep=comma]{NSE_WSODIRK433.csv};
          \addplot[pink, mark=o, mark size = 3, mark options = {solid}, dashed, line width = 2]
          table [x=dt,y=elt, col sep=comma]{NSE_WSODIRK744.csv};
          \addplot[domain=0.002:0.02] {120/pow(x,-5)} node[below, midway, yshift=-1pt, anchor=west]{$\mathcal{O}(\Delta t^5)$};
          \addplot[domain=0.002:0.02] {1000/pow(x,-3)} node[below, midway, yshift=2pt, anchor=east]{$\mathcal{O}(\Delta t^3)$};
        \end{axis}
      \end{tikzpicture}
      \caption{Error in final lift}
      \label{fig:nseelt}           
  \end{subfigure}    
  \end{center} 
  \caption{Error in drag and lift coefficients as a function of time step size for various Runge--Kutta schemes.}
  \label{fig:errnse}
\end{figure}
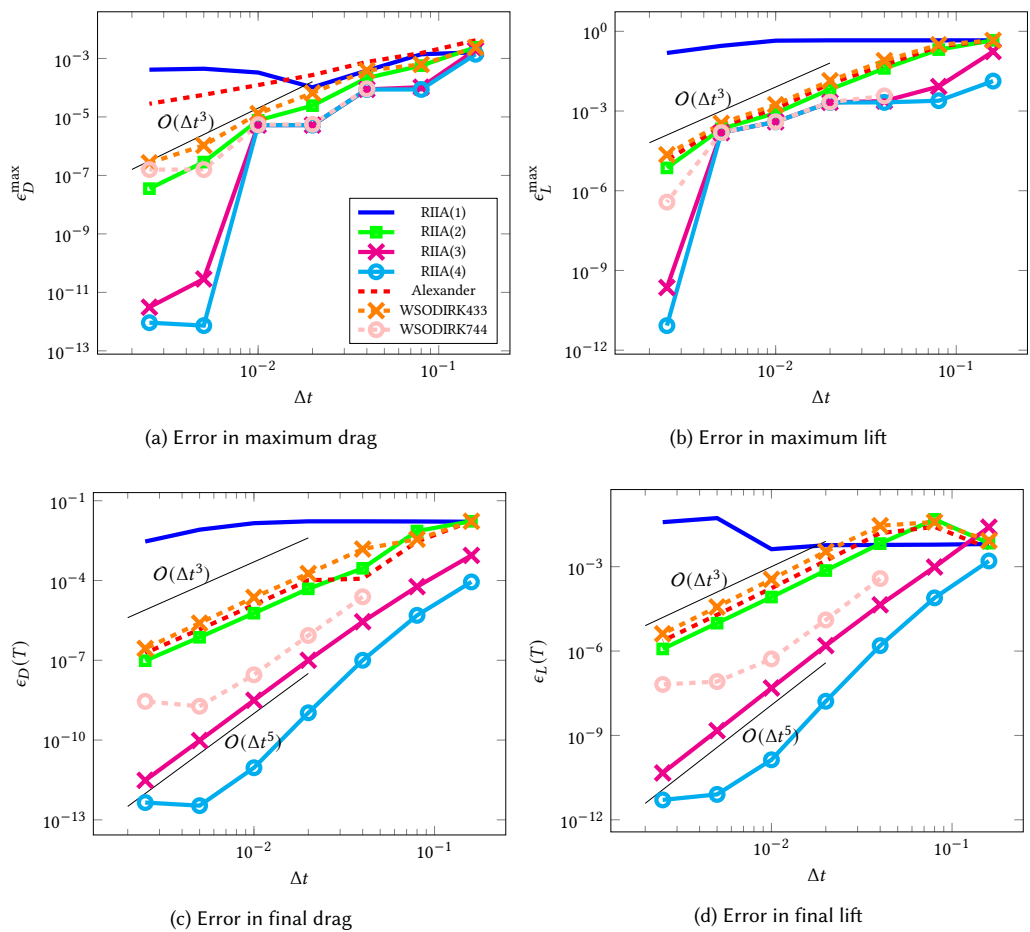

\begin{figure}
  \begin{center}
    \begin{subfigure}[c]{0.45\textwidth}
      \begin{tikzpicture}[scale=0.8]
        \begin{axis}[xlabel=Time(s), ylabel={$\epsilon_D^{\max}$}, xmode=log,ymode=log, legend pos=south west, legend style={nodes={scale=0.8}}]
          \addplot[blue, line width = 2]
          table [x=time,y=edmax, col sep=comma]{NSE_RadauIIA1.csv};
          \addlegendentry{RIIA(1)}
          \addplot[green, mark=square, mark options = {solid}, line width = 2]
          table [x=time,y=edmax, col sep=comma]{NSE_RadauIIA2.csv};
          \addlegendentry{RIIA(2)}
          \addplot[magenta, mark=x, mark size = 5, mark options = {solid}, line width = 2]
          table [x=time,y=edmax, col sep=comma]{NSE_RadauIIA3.csv};
          \addlegendentry{RIIA(3)}
          \addplot[cyan, mark=o, mark size = 3, mark options = {solid}, line width = 2]
          table [x=time,y=edmax, col sep=comma]{NSE_RadauIIA4.csv};
          \addlegendentry{RIIA(4)}
          \addplot[red, dashed, line width = 2]
          table [x=time,y=edmax, col sep=comma]{NSE_Alexander.csv};
          \addlegendentry{Alexander}
          \addplot[orange, mark=x, mark size = 5, mark options = {solid}, dashed, line width = 2]
          table [x=time,y=edmax, col sep=comma]{NSE_WSODIRK433.csv};
          \addlegendentry{WSODIRK433}
          \addplot[pink, mark=o, mark size = 3, mark options = {solid}, dashed, line width = 2]
          table [x=time,y=edmax, col sep=comma]{NSE_WSODIRK744.csv};
          \addlegendentry{WSODIRK744}
        \end{axis}
      \end{tikzpicture}
      \caption{Error in maximum drag}
      \label{fig:nseedmaxvtime}
    \end{subfigure}
    \begin{subfigure}[c]{0.45\textwidth}
      \begin{tikzpicture}[scale=0.8]
        \begin{axis}[xlabel=Time (s), ylabel={$\epsilon_L^{\max}$},xmode=log, ymode=log]
        \addplot[blue, line width = 2]
          table [x=time,y=elmax, col sep=comma]{NSE_RadauIIA1.csv};
          \addplot[green, mark=square, mark options = {solid}, line width = 2]
          table [x=time,y=elmax, col sep=comma]{NSE_RadauIIA2.csv};
          \addplot[magenta, mark=x, mark size = 5, mark options = {solid}, line width = 2]
          table [x=time,y=elmax, col sep=comma]{NSE_RadauIIA3.csv};
          \addplot[cyan, mark=o, mark size = 3, mark options = {solid}, line width = 2]
          table [x=time,y=elmax, col sep=comma]{NSE_RadauIIA4.csv};
          \addplot[red, dashed, line width = 2]
          table [x=time,y=elmax, col sep=comma]{NSE_Alexander.csv};
          \addplot[orange, mark=x, mark size = 5, mark options = {solid}, dashed, line width = 2]
          table [x=time,y=elmax, col sep=comma]{NSE_WSODIRK433.csv};
          \addplot[pink, mark=o, mark size = 3, mark options = {solid}, dashed, line width = 2]
          table [x=time,y=elmax, col sep=comma]{NSE_WSODIRK744.csv};
        \end{axis}
      \end{tikzpicture}
      \caption{Error in maximum lift}
      \label{fig:nseelmaxvtime}      
    \end{subfigure} \\[12pt]
    \begin{subfigure}[c]{0.45\textwidth}
      \begin{tikzpicture}[scale=0.8]
        \begin{axis}[xlabel=Time (s), ylabel={$\epsilon_D(T)$}, xmode=log,ymode=log]
          \addplot[blue, line width = 2]
          table [x=time,y=edt, col sep=comma]{NSE_RadauIIA1.csv};
          \addplot[green, mark=square, mark options = {solid}, line width = 2]
          table [x=time,y=edt, col sep=comma]{NSE_RadauIIA2.csv};
          \addplot[magenta, mark=x, mark size = 5, mark options = {solid}, line width = 2]
          table [x=time,y=edt, col sep=comma]{NSE_RadauIIA3.csv};
          \addplot[cyan, mark=o, mark size = 3, mark options = {solid}, line width = 2]
          table [x=time,y=edt, col sep=comma]{NSE_RadauIIA4.csv};
          \addplot[red, dashed, line width = 2]
          table [x=time,y=edt, col sep=comma]{NSE_Alexander.csv};
          \addplot[orange, mark=x, mark size = 5, mark options = {solid}, dashed, line width = 2]
          table [x=time,y=edt, col sep=comma]{NSE_WSODIRK433.csv};
          \addplot[pink, mark=o, mark size = 3, mark options = {solid}, dashed, line width = 2]
          table [x=time,y=edt, col sep=comma]{NSE_WSODIRK744.csv};
        \end{axis}
      \end{tikzpicture}
      \caption{Error in final drag}
      \label{fig:nseedtvtime}     
    \end{subfigure}
    \begin{subfigure}[c]{0.45\textwidth}
      \begin{tikzpicture}[scale=0.8]
        \begin{axis}[xlabel=Time(s), ylabel={$\epsilon_L(T)$},xmode=log, ymode=log, legend pos=south west]
        \addplot[blue, line width = 2]
          table [x=time,y=elt, col sep=comma]{NSE_RadauIIA1.csv};
          \addplot[green, mark=square, mark options = {solid}, line width = 2]
          table [x=time,y=elt, col sep=comma]{NSE_RadauIIA2.csv};
          \addplot[magenta, mark=x, mark size = 5, mark options = {solid}, line width = 2]
          table [x=time,y=elt, col sep=comma]{NSE_RadauIIA3.csv};
          \addplot[cyan, mark=o, mark size = 3, mark options = {solid}, line width = 2]
          table [x=time,y=elt, col sep=comma]{NSE_RadauIIA4.csv};
          \addplot[red, dashed, line width = 2]
          table [x=time,y=elt, col sep=comma]{NSE_Alexander.csv};
          \addplot[orange, mark=x, mark size = 5, mark options = {solid}, dashed, line width = 2]
          table [x=time,y=elt, col sep=comma]{NSE_WSODIRK433.csv};
          \addplot[pink, mark=o, mark size = 3, mark options = {solid}, dashed, line width = 2]
          table [x=time,y=elt, col sep=comma]{NSE_WSODIRK744.csv};
        \end{axis}
      \end{tikzpicture}
      \caption{Error in final lift}
      \label{fig:nseeltvtime}           
  \end{subfigure}    
  \end{center} 
  \caption{Error in drag and lift coefficients as a function of run time for various Runge--Kutta schemes.}
  \label{fig:errnsevtime}
\end{figure}
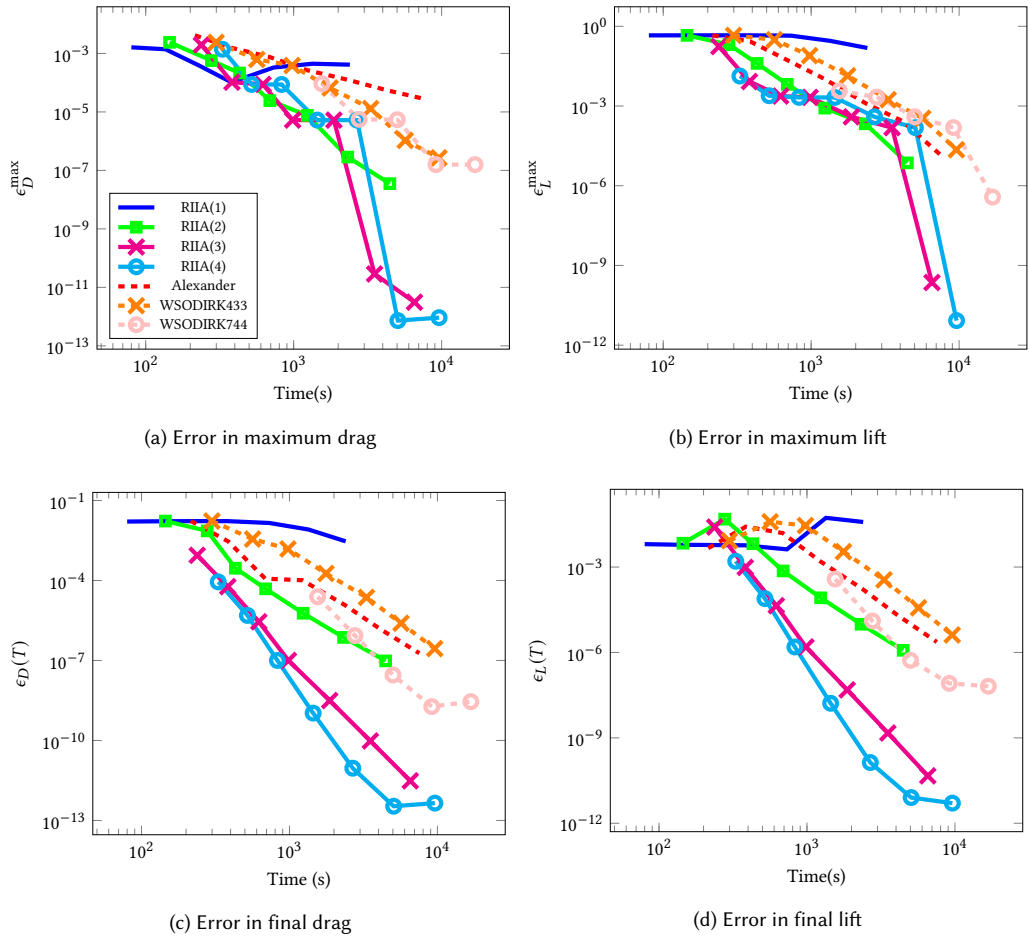

The overall solver performance is reported in Figure~\ref{fig:perfnse}, where we see several features.
First, for a fixed number of stages, we see that both the number of Newton iterations per time step and the number of Krylov iterations per Newton step decrease slightly with $\Delta t$, which is perhaps expected.
Fixing $\Delta t$ and varying the number of stages, we see that the Newton and Krylov iteration counts remain very stable.
This is consistent with results reported in~\cite{abu2022monolithic, mmg}.
Second, for a fixed number of stages, the run time is slightly less than proportional to the number of time steps taken (reciprocal of $\Delta t$) owing to the improvements in linear and nonlinear convergence.  
Given the convergence behavior, a run-time roughly proportional to the number of stages is expected for DIRKs.
This pattern is also observed for the monolithic multigrid solver for RadauIIA methods.
The three-stage method for a given time step takes slightly more than three times the single stage method, and the four-stage method takes slightly less than four times.
It is also important to note that DIRKs with a given number of stages give run-time very close to RadauIIA methods with the same number of stages.
So, although the fully implicit method may use more memory than DIRKs, they do not seem to lead to worse run-times.

It is equally important to observe the accuracy of our methods in various metrics, as reported in Figure~\ref{fig:errnse}.
Fixing the mesh and reducing the time step means that the effective stiffness of the system is reduced for small $\Delta t$, and so we may expect to see order reduction decrease in this limit.
We notice barely any convergence of the single-stage method as $\Delta t$ decreases.
Perhaps first-order asymptotic convergence would be observed with further reduction in $\Delta t$.  We also note that the RadauIIA(3) and RadauIIA(4) methods produce essentially identical values for maximum drag and lift as our reference solution at the finest time step (or two).  This is likely due to the fact that our reference solution is computed with RadauIIA(5) at the same (finest) time step, but we are unable to produce an independent reference value to improve accuracy of these measurements.

We have three methods with formal accuracy of order three -- RadauIIA(2), Alexander, and WSODIRK433.
Among these, comparable accuracy is obtained, with the RadauIIA(2) method slightly better than Alexander, which is in turn slightly better than WSODIRK433.  WSODIRK744 has formal accuracy of order 4, while RadauIIA(3) has formal accuracy of order 5.  Among these, RadauIIA(3) has a slight advantage in accuracy.
Both of these methods give slightly worse accuracy compared to RadauIIA(4) at the final time step and similar accuracy for the maximum lift/drag values.  We note, however, that for an index 2 DAE, like Navier-Stokes, the accuracy guarantees are given by the stage order of the schemes, and that DIRK schemes are only guaranteed to have stage order one or two, while RadauIIA($s$) has stage order $s$.  In general, we see the schemes outperform the stage-order bound for the error in final drag and lift, achieving their formal order of accuracy.  Results are slightly less clear for the maximum drag and lift, although we clearly see that Alexander under-achieves its formal order of accuracy for the maximum drag.

Our computations also let us compare accuracy achieved versus overall run-time, which is reported in Figure~\ref{fig:errnsevtime}.
For lift/drag error at the final time, the RadauIIA methods with 3 and 4 stages clearly produce the best results.  Although WSODIRK744 is quite accurate, it is far more expensive to compute.
We also note that using the Rana preconditioner rather than monolithic multigrid scheme for the 2- and 3- stage Radau methods would further tilt the accuracy versus run-time consideration in the favor of fully implicit schemes.

\subsection{Cahn--Hilliard}
As a further example of our advances in \Irksome{}, we consider the Cahn--Hilliard equation, which models phase separation in a binary fluid and is given by
\begin{equation}
  \label{eq:ch}
  c_t - \nabla \cdot M \left( \nabla \left( f^\prime(c) - \lambda \Delta c \right) \right)  = 0.
\end{equation}
The dependent variable, $c$, controls the relative mixing of the two phases of the fluid, with $c=\pm 1$ indicating a single phase is present and values in between indicating some degree of mixture.
The scalar $M$ is the \emph{mobility} and can, in general, depend on $c$ (although, for simplicity, we do not consider that case here), and $\lambda$ is a constant determining the size of free energy for a given concentration gradient.  The function $f$ is generally taken to be a non-convex function (e.g. double well potential).  In our case, we take the qualitative choice of
\begin{equation}
  \label{eq:fc}
  f(c) = 50 (c^2 - 1)^2.
\end{equation}
We close the system with boundary conditions
\begin{equation}
  \label{eq:chbc}
  \begin{split}
M\left(\nabla\left(f^\prime(c) - \lambda \Delta c\right)\right) \cdot n &= 0 \quad {\rm on} \ \partial\Omega, \\
M \lambda \nabla c \cdot n &= 0 \quad {\rm on} \ \partial\Omega.
  \end{split}
\end{equation}
A primary physical quantity for the Cahn--Hilliard problem is the free energy
\begin{equation}
  \label{eq:freeenergy}
  E(c) = \int_\Omega M f(c) + \tfrac{M \lambda}{2} |\nabla c|^2 \, dx,
\end{equation}
which is nonincreasing in time, and it is highly desirable for numerical methods to preserve this feature.  In the example below we take $M=1$, $\lambda = 10^{-2}$.

The fourth-order derivatives in~\eqref{eq:ch} require special care -- either using high-continuity ($C^1$) elements or an interior penalty technique~\cite{wells2006discontinuous}, or writing the problem as a system of second-order equations to use standard $C^0$ elements.
While the latter two options are easily realized, Firedrake also supports certain $C^1$ elements on triangles~\cite{finat-zany}.  Here, we present some two-dimensional examples using a conforming primal method and focus on issues related to time-stepping and energy stability.
Thus, we consider an $H^2$-conforming finite element space, $V_h$, and write the weak form of~\eqref{eq:ch} to find $c:[0, T] \rightarrow V_h$ such that, for almost all $t$,
\begin{equation}
  \label{eq:chweak}
  \left( c_t, v \right) +
  M \left( \nabla f^\prime(c), \nabla v \right))
  + M \lambda \left( \Delta c, \Delta v \right)
  - M \lambda \langle \Delta c, \tfrac{\partial v}{\partial n} \rangle
  - M \lambda \langle \tfrac{\partial c}{\partial n} , \Delta v \rangle
  + \tfrac{\beta M \lambda}{h} \langle \tfrac{\partial c}{\partial n}, \tfrac{\partial v}{\partial n} \rangle
  = 0
\end{equation}
for all $v \in V_h$.  Here, $\langle \cdot, \cdot \rangle$ indicates the $L^2$ inner product over the boundary of $\Omega$, while $( \cdot,\cdot )$ continues to denote the $L^2$ inner product over all of $\Omega$.

As in~\cite{finat-zany}, we use the Bell element~\cite{bell1969refined}, an $H^2$ conforming triangle consisting of quintic polynomials with normal derivatives on each edge restricted to univariate cubics.  The second condition in~\eqref{eq:chbc} cannot be directly implemented strongly, so we augment the variational form with a Nitsche-type technique, choosing the parameter $\beta$ to be sufficiently large; here, $\beta = 250$.
However, the first boundary condition in~\eqref{eq:chbc} arises via integration by parts as a natural boundary condition.

We consider the Cahn--Hilliard equation on a unit square domain, divided into
a 32 $\times$ 32 mesh of squares subdivided into right triangles, setting the initial condition to 0.42 on $[0.3,0.7]^2$ and 0.38 on the balance of the unit square.
We evolve the semidiscrete problem~\eqref{eq:chweak} forward in time using RadauIIA($s$) methods with $1 \leq s \leq 4$ as well as with the Alexander, WSODIRK433, and WSODIRK744 methods and various time step sizes.
At each time step, we perform a Newton iteration with $10^{-10}$ stopping tolerances on both the absolute and relative values of the nonlinear residual, solving the Jacobian system with FGMRES preconditioned with stage-coupled monolithic geometric multigrid.
As relaxation, we use a vertex patch-based additive Schwarz method, applied with three Chebyshev iterations, estimating the largest eigenvalue, $\mu$, of the additive Schwarz preconditioned system and using Chebyshev polynomials defined on the interval $[0.25\mu,1.1\mu]$. 
On the coarsest grid, we use a sparse direct method.  The FGMRES tolerance is chosen adaptively according to the Eisenstat--Walker technique~\cite{eisenstat1996choosing}.
We applied this multigrid method both monolithically and for each stage in a Rana-type preconditioner.  We found that for smaller time steps, the monolithic scheme gave better run-time, although the reverse was true for larger time steps.  Here, we report only on monolithic methods.

\begin{figure}
  \begin{subfigure}[c]{0.475\textwidth}
    \includegraphics[width=0.9\textwidth]{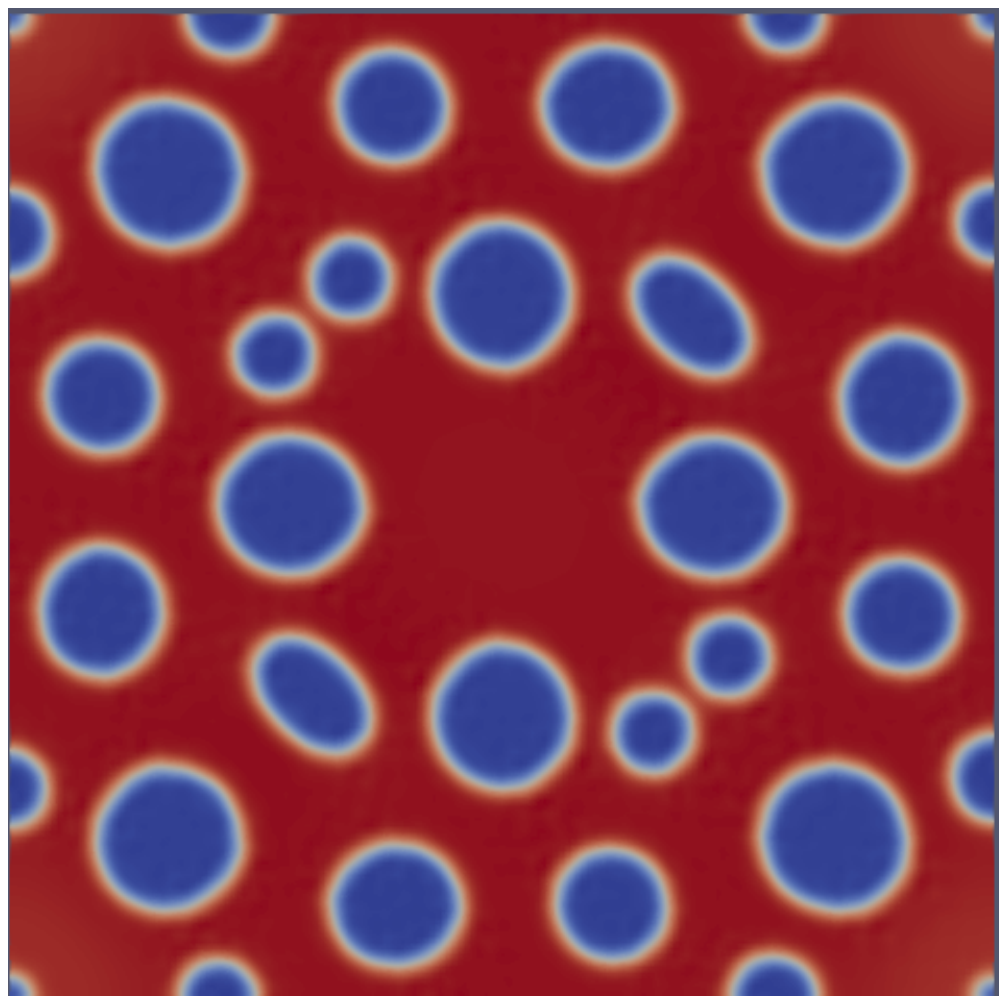}
    \caption{RadauIIA(2)}
  \end{subfigure}
  \begin{subfigure}[c]{0.475\textwidth}
    \includegraphics[width=0.9\textwidth]{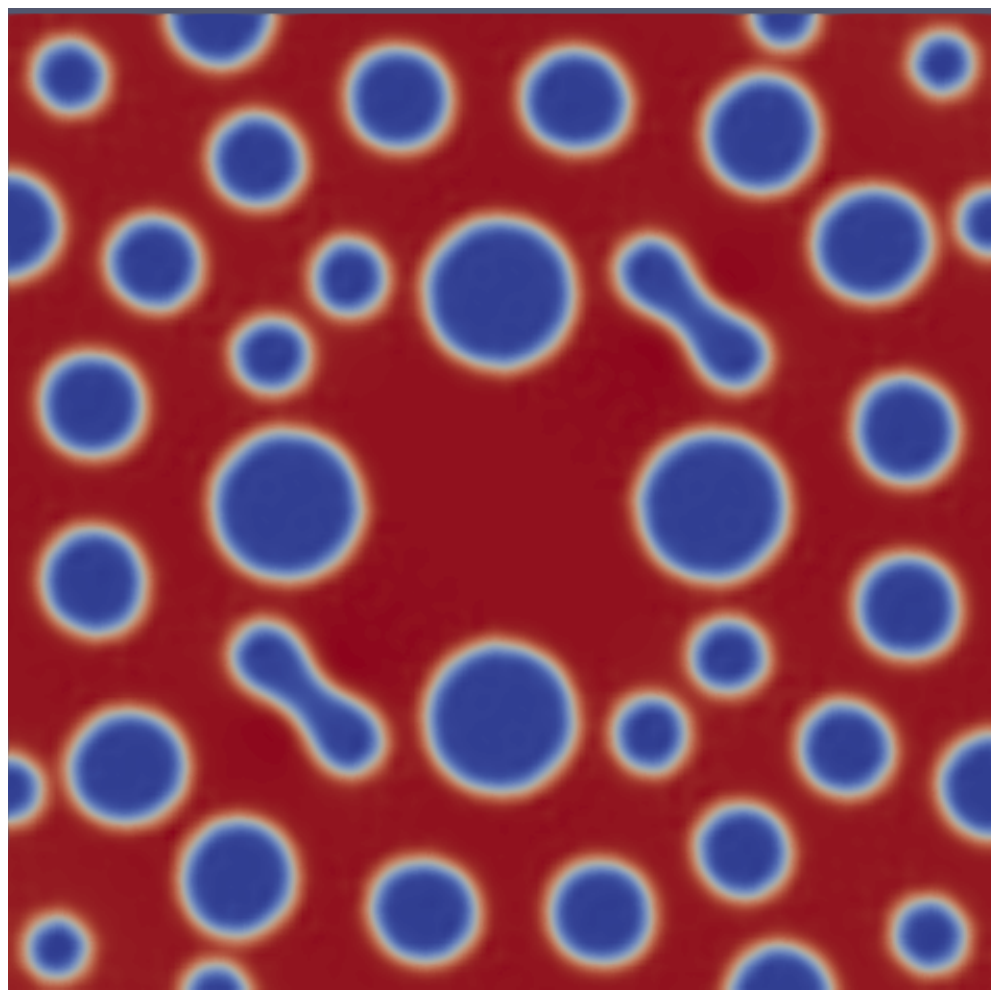}
    \caption{RadauIIA(3)}
  \end{subfigure}
  \caption{Final Cahn-Hilliard states computed with 2- and 3-stage RadauIIA time integration.}
  \label{fig:chsol}
\end{figure}

Our results suggest that increasing the number of stages can have a positive effect on both accuracy and efficiency in our methods.
In Figure~\ref{fig:chsol}, we show the results of RadauIIA(2) and RadauIIA(3), both with $\Delta t= 2.5\times 10^{-6}$.  Refining $\Delta t$ in the 2-stage method produces a final result quite similar to RadauIIA(3), while reducing the time step for the 3-stage method or using a 4-stage method does not change the result significantly.

\begin{figure}
    \begin{center}
  \begin{tikzpicture}[scale=0.85]
    \begin{axis}[xlabel={$t$}, ylabel={$F(c)$}, ymin=0, ymax=40,
        legend pos = outer north east]
      \addplot[blue, line width = 2]
      table [x=t,y=E, col sep=comma]{ch.energy.RIIA1.IA.6.250e-07.csv};
      \addlegendentry{RIIA(1), $\Delta t$=6.25e-07}
      \addplot[green, mark=square, mark options = {solid}, line width = 2, mark repeat = 20]
      table [x=t,y=E, col sep=comma]{ch.energy.RIIA2.IA.1.250e-06.csv};
      \addlegendentry{RIIA(2), $\Delta t$=1.25e-6}
      \addplot[magenta, mark=x, mark size = 5, mark options = {solid}, line width = 2, mark repeat = 40, mark phase = 5]
      table [x=t,y=E, col sep=comma]{ch.energy.RIIA3.IA.2.500e-06.csv};
      \addlegendentry{RIIA(3), $\Delta t$=2.5e-6}
      \addplot[cyan, mark=o, mark size = 3, mark options = {solid}, line width = 2, mark repeat = 25]
      table [x=t,y=E, col sep=comma]{ch.energy.RIIA4.IA.5.000e-06.csv};
      \addlegendentry{RIIA(4), $\Delta t$=5.0e-6}
      \addplot[red, dashed, line width = 2]
      table [x=t,y=E, col sep=comma]{ch.energy.Alexander.1.250e-06.csv};
      \addlegendentry{Alexander, $\Delta t$=1.25e-6}
      \addplot[orange, mark=x, mark size = 5, mark options = {solid}, dashed, line width = 2, mark repeat = 40, mark phase = 10]
      table [x=t,y=E, col sep=comma]{ch.energy.WSODIRK433.1.250e-06.csv};
      \addlegendentry{WSODIRK433, $\Delta t$=1.25e-6}
      \addplot[pink, mark=o, mark size = 3, mark options = {solid}, dashed, line width = 2, mark repeat = 40, mark phase = 20]
      table [x=t,y=E, col sep=comma]{ch.energy.WSODIRK744.1.250e-06.csv};
      \addlegendentry{WSODIRK744, $\Delta t$=1.25e-6}      
    \end{axis}
  \end{tikzpicture}
  \end{center}
  \caption{Ginzburg-Landau free energy versus time for the Cahn-Hilliard equation for various numbers of stages/time steps.}
  \label{fig:energy}
\end{figure}
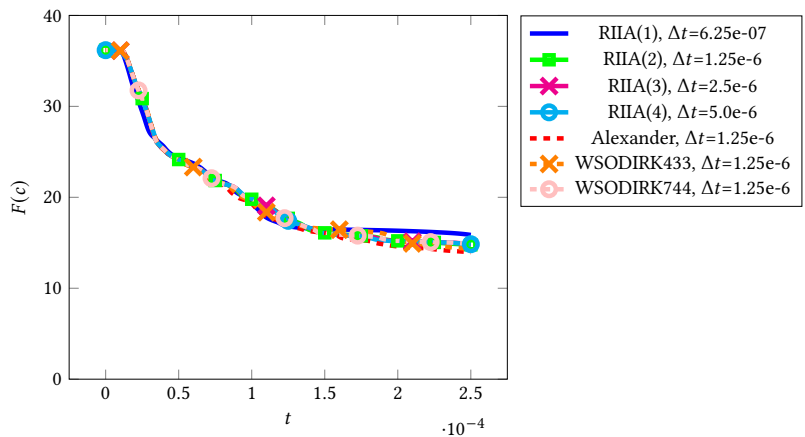

We also consider the energy from~\eqref{eq:freeenergy}.
Figure~\ref{fig:energy} shows a monotonically decreasing free energy for each method, but we see that the 1-stage method (implicit Euler), even with $\Delta t = 6.25 \times 10^{-7}$, gives a qualitatively different result compared to the other cases.  We see increased accuracy with a smaller $\Delta t$, in line with the results provided with three- and four-stage methods with slightly larger times.

We collected statistics on the solver -- overall run-time and number of linear and nonlinear iterations as well as the relative error in the final energy denoted by $\epsilon(T)$.
We used time steps from $\Delta t=5\times 10^{-6}$ down to $3.125 \times 10^{-7}$.
As a reference value, we used the energy computed with a RadauIIA(5) method with time step of $3.125 \times 10^{-7}$ .
Some of the methods experienced failure in Newton convergence for the larger time steps considered, so some lines in the plots have fewer data points than others.

\begin{figure}
  \begin{center}
    \begin{subfigure}[c]{0.3\textwidth}
    \begin{tikzpicture}[scale=0.5]
      \begin{semilogxaxis}[xlabel=$\Delta t$, ylabel={Time (s)},
          ymode=log,
            label style={font=\LARGE},
            tick label style={font=\LARGE}]
        \addplot[blue, line width = 2]
        table [x=dt,y=time, col sep=comma]{ch.RIIA1.csv};
        \addplot[green, mark=square, mark options = {solid}, line width = 2]
        table [x=dt,y=time, col sep=comma]{ch.RIIA2.csv};
        \addplot[magenta, mark=x, mark size = 5, mark options = {solid}, line width = 2]
        table [x=dt,y=time, col sep=comma]{ch.RIIA3.csv};
        \addplot[cyan, mark=o, mark size = 3, mark options = {solid}, line width = 2]
        table [x=dt,y=time, col sep=comma]{ch.RIIA4.csv};
        \addplot[red, dashed, line width = 2]
        table [x=dt,y=time, col sep=comma]{ch.Alexander.csv};
        \addplot[orange, mark=x, mark size = 5, mark options = {solid}, dashed, line width = 2]
        table [x=dt,y=time, col sep=comma]{ch.WSODIRK433.csv};
        \addplot[pink, mark=o, mark size = 3, mark options = {solid}, dashed, line width = 2]
        table [x=dt,y=time, col sep=comma]{ch.WSODIRK744.csv};
      \end{semilogxaxis}
    \end{tikzpicture}
    \caption{Timing}
    \label{fig:chtime}
    \end{subfigure}
    \begin{subfigure}[c]{0.3\textwidth}
      \begin{tikzpicture}[scale=0.5]
        \begin{axis}[xlabel=$\Delta t$, ylabel={Nonlinear Iterations},ymin=0, xmode=log,
            label style={font=\LARGE},
            tick label style={font=\LARGE}]
          \addplot[blue, line width = 2]
          table [x=dt,y=newt, col sep=comma]{ch.RIIA1.csv};
          \addplot[green, mark=square, mark options = {solid}, line width = 2]
          table [x=dt,y=newt, col sep=comma]{ch.RIIA2.csv};
          \addplot[magenta, mark=x, mark size = 5, mark options = {solid}, line width = 2]
          table [x=dt,y=newt, col sep=comma]{ch.RIIA3.csv};
          \addplot[cyan, mark=o, mark size = 3, mark options = {solid}, line width = 2]
          table [x=dt,y=newt, col sep=comma]{ch.RIIA4.csv};
          \addplot[red, dashed, line width = 2]
          table [x=dt,y=newt, col sep=comma]{ch.Alexander.csv};
          \addplot[orange, mark=x, mark size = 5, mark options = {solid}, dashed, line width = 2]
          table [x=dt,y=newt, col sep=comma]{ch.WSODIRK433.csv};
          \addplot[pink, mark=o, mark size = 3, mark options = {solid}, dashed, line width = 2]
          table [x=dt,y=newt, col sep=comma]{ch.WSODIRK744.csv};
        \end{axis}
      \end{tikzpicture}
      \caption{Newton iterations}
      \label{fig:chnewt}
    \end{subfigure}
  \begin{subfigure}[c]{0.3\textwidth}
    \begin{tikzpicture}[scale=0.5]
      \begin{axis}[xlabel=$\Delta t$, ylabel={GMRES Iterations}, ymin=0,xmode=log, legend style={nodes={scale=0.8},font=\LARGE},
        legend pos = outer north east,
            label style={font=\LARGE},
            tick label style={font=\LARGE}]
        \addplot[blue, line width = 2]
          table [x=dt,y=ksp, col sep=comma]{ch.RIIA1.csv};
          \addlegendentry{RIIA(1)}
          \addplot[green, mark=square, mark options = {solid}, line width = 2]
          table [x=dt,y=ksp, col sep=comma]{ch.RIIA2.csv};
          \addlegendentry{RIIA(2)}
          \addplot[magenta, mark=x, mark size = 5, mark options = {solid}, line width = 2]
          table [x=dt,y=ksp, col sep=comma]{ch.RIIA3.csv};
          \addlegendentry{RIIA(3)}
          \addplot[cyan, mark=o, mark size = 3, mark options = {solid}, line width = 2]
          table [x=dt,y=ksp, col sep=comma]{ch.RIIA4.csv};
          \addlegendentry{RIIA(4)}
          \addplot[red, dashed, line width = 2]
          table [x=dt,y=ksp, col sep=comma]{ch.Alexander.csv};
          \addlegendentry{Alexander}
          \addplot[orange, mark=x, mark size = 5, mark options = {solid}, dashed, line width = 2]
          table [x=dt,y=ksp, col sep=comma]{ch.WSODIRK433.csv};
          \addlegendentry{WSODIRK433}
          \addplot[pink, mark=o, mark size = 3, mark options = {solid}, dashed, line width = 2]
          table [x=dt,y=ksp, col sep=comma]{ch.WSODIRK744.csv};
          \addlegendentry{WSODIRK744}
        \end{axis}
    \end{tikzpicture}
    \caption{GMRES iterations}
    \label{fig:chgmres}
  \end{subfigure}
  \end{center}
  \caption{Solver performance for integrating Cahn-Hilliard equations on $[0,T]$.  We report total wall-clock time for the time stepping loop as well as average nonlinear iterations per time step and average linear iterations per nonlinear solve.  For DIRKs, our Newton iterations are averaged over each nonlinear stage of each time step.}
  \label{fig:perfch}
\end{figure}
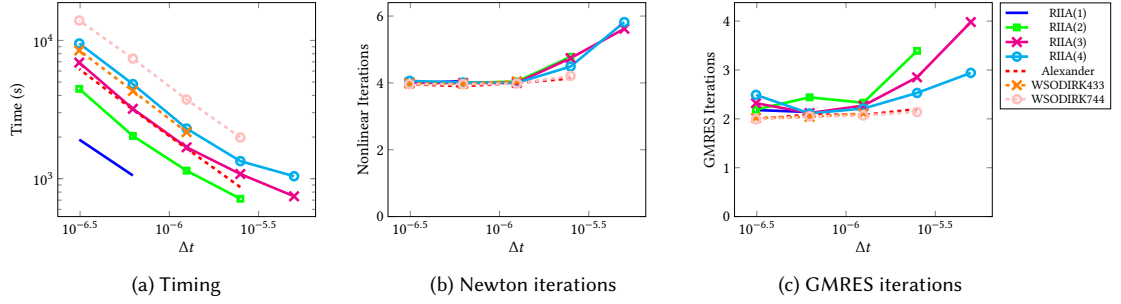

In Figure~\ref{fig:chtime}, we see the run time of each method versus time step size.
These timings show that run-time depends primarily on the number of stages.
Unlike the heat equation and Navier--Stokes, here we note that RadauIIA requires slightly more run-time than a DIRK with the same number of stages.
Figures~\ref{fig:chnewt} and~\ref{fig:chgmres} indicate the number of Newton and GMRES iterations required for each method for various time step sizes.  As the time step decreases, we approach about four nonlinear iterations per time step (per stage as well for DIRKs), and the number of linear iterations per nonlinear solve is very small.

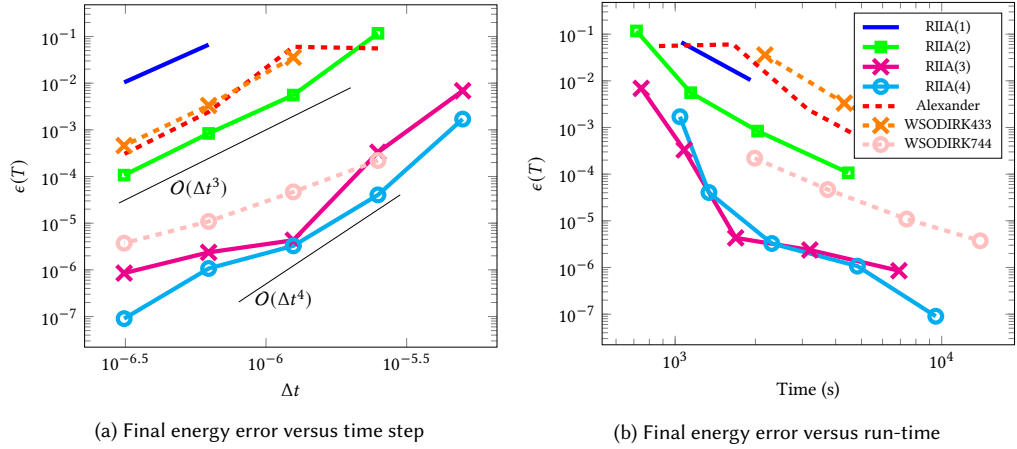
\begin{figure}
  \begin{center}
    \begin{subfigure}[c]{0.45\textwidth}
      \begin{tikzpicture}[scale=0.8]
        \begin{axis}[xlabel=$\Delta t$, ylabel={$\epsilon(T)$}, xmode=log, ymode=log]
          \addplot[blue, line width = 2]
          table [x=dt,y=et, col sep=comma]{ch.RIIA1.csv};
          \addplot[green, mark=square, mark options = {solid}, line width = 2]
          table [x=dt,y=et, col sep=comma]{ch.RIIA2.csv};
          \addplot[magenta, mark=x, mark size = 5, mark options = {solid}, line width = 2]
          table [x=dt,y=et, col sep=comma]{ch.RIIA3.csv};
          \addplot[cyan, mark=o, mark size = 3, mark options = {solid}, line width = 2]
          table [x=dt,y=et, col sep=comma]{ch.RIIA4.csv};
          \addplot[red, dashed, line width = 2]
          table [x=dt,y=et, col sep=comma]{ch.Alexander.csv};
          \addplot[orange, mark=x, mark size = 5, mark options = {solid}, dashed, line width = 2]
          table [x=dt,y=et, col sep=comma]{ch.WSODIRK433.csv};
          \addplot[pink, mark=o, mark size = 3, mark options = {solid}, dashed, line width = 2]
          table [x=dt,y=et, col sep=comma]{ch.WSODIRK744.csv};
          \addplot[domain=0.0000003:0.000002] {1000000000000000/pow(x,-3)} node[below, midway, yshift=-20pt, anchor=east]{$\mathcal{O}(\Delta t^3)$};
          \addplot[domain=0.0000008:0.000003] {500000000000000000/pow(x,-4)} node[below, midway, yshift=-24pt, anchor=east]{$\mathcal{O}(\Delta t^4)$};
        \end{axis}
      \end{tikzpicture}
      \caption{Final energy error versus time step}
      \label{fig:chetvdt}
    \end{subfigure}
    \begin{subfigure}[c]{0.45\textwidth}
      \begin{tikzpicture}[scale=0.8]
        \begin{axis}[xlabel=Time (s), ylabel={$\epsilon(T)$}, xmode=log, ymode=log, legend pos=north east, legend style={nodes={scale=0.8}}]
          \addplot[blue, line width = 2]
          table [x=time,y=et, col sep=comma]{ch.RIIA1.csv};
          \addlegendentry{RIIA(1)}
          \addplot[green, mark=square, mark options = {solid}, line width = 2]
          table [x=time,y=et, col sep=comma]{ch.RIIA2.csv};
          \addlegendentry{RIIA(2)}
          \addplot[magenta, mark=x, mark size = 5, mark options = {solid}, line width = 2]
          table [x=time,y=et, col sep=comma]{ch.RIIA3.csv};
          \addlegendentry{RIIA(3)}
          \addplot[cyan, mark=o, mark size = 3, mark options = {solid}, line width = 2]
          table [x=time,y=et, col sep=comma]{ch.RIIA4.csv};
          \addlegendentry{RIIA(4)}
          \addplot[red, dashed, line width = 2]
          table [x=time,y=et, col sep=comma]{ch.Alexander.csv};
          \addlegendentry{Alexander}
          \addplot[orange, mark=x, mark size = 5, mark options = {solid}, dashed, line width = 2]
          table [x=time,y=et, col sep=comma]{ch.WSODIRK433.csv};
          \addlegendentry{WSODIRK433}
          \addplot[pink, mark=o, mark size = 3, mark options = {solid}, dashed, line width = 2]
          table [x=time,y=et, col sep=comma]{ch.WSODIRK744.csv};
          \addlegendentry{WSODIRK744}
        \end{axis}
      \end{tikzpicture}
      \caption{Final energy error versus run-time}
      \label{fig:chetvtime}
    \end{subfigure}    
  \end{center}
  \caption{Error in final free-energy.}
  \label{fig:cherr}
\end{figure}

It is also interesting to measure the accuracy obtained, both in absolute terms and relative to computational effort.
Figure~\ref{fig:cherr} reports on these.
In Figure~\ref{fig:chetvdt} we see that, in absolute terms, the RadauIIA methods with three and four stages give the lowest error over most time steps considered, followed by WSODIRK744.
Here, RadauIIA(2) gives better accuracy than either the Alexander or WSODIRK433 methods (and at much lower cost).
Figure~\ref{fig:chetvtime} plots the accuracy obtained versus the run-time for the simulation.  Here, we see that the three- and four-stage RadauIIA methods perform very well.
The WSODIRK744 method comes in next -- although it is the most expensive, it is also very accurate.
RadauIIA(2) beats both of the other DIRKs considered.  Again we note the difference between formal order of accuracy and that obtained by the stage order, due to stiffness of the resulting method-of-lines discretization, so that few of the schemes appear to attain their formal order of accuracy.

\section{Conclusions}\label{sec:conclusions}
In this work, we demonstrate several important improvements in the \Irksome{} package developed since publication of~\cite{farrell2021irksome}.  This includes higher-accuracy treatment of strong Dirichlet boundary conditions, alternative stage-derivative and stage-value formulations of Runge-Kutta methods, better support for DIRK schemes, and enhancements in expressing effective preconditioners through Firedrake and PETSc options.  Taken together, these are demonstrated to lead to a much more capable set of tools to enable higher accuracy and solution of more complex systems of differential and differential-algebraic equations.  Numerical results in this paper show the efficiency of some of the improved preconditioning options, as well as the accuracy that can be obtained using the fully implicit Runge-Kutta methodology over comparable DIRK schemes.

A focus for future work is on continuing to grow the support for higher fidelity numerical methods and higher efficiency solvers.
One key effort currently underway is the development of a reliable adaptive time-stepping scheme for implicit RK discretizations of systems of DAEs, and its implementation into the \Irksome{} package.
Future work also includes the development of similar support for IMEX methods, where part of the system is treated implicitly and part explicitly, and evaluating their performance relative to fully implicit methods.

\appendix

\section{Some Butcher tableaux} \label{sec:butcher}
For reference, we have included several Butcher tableaux of interest.
The RadauIIA($s$) methods are all continuous collocation schemes, with the $c_i$ values taken as roots of
\begin{equation}
  \frac{d^{s-1}}{dx^{s-1}} \left( x^{s-1} (x-1)^s \right).
\end{equation}
The values of $b_i$ and $A_{ij}$ follow programmatically from this choice of points.  For this construction, we refer the reader to~\cite{guillou1969resolution, haireri}. 

The RadauIIA(1) method, which coincides with backward Euler, is given by the tableaux
\begin{equation}
  \begin{array}{c|c}
    1 & 1 \\ \hline
    & 1
  \end{array}
\end{equation}
The two-stage RadauIIA method has the tableau
\begin{equation}
  \begin{array}{c|cc}
    \tfrac{1}{3} & \tfrac{5}{12} & -\tfrac{1}{12} \\
    1 & \tfrac{3}{4} & \tfrac{1}{4} \\ \hline
    & \tfrac{3}{4} & \tfrac{1}{4} \\
  \end{array},
\end{equation}
    and that for the three-stage method is
\begin{equation}
  \begin{array}{c|ccc}
    \tfrac{2}{5} - \tfrac{\sqrt{6}}{10} & \tfrac{11}{45} - \tfrac{7\sqrt{6}}{360} &
    \tfrac{37}{225} - \tfrac{169\sqrt{6}}{1800} & -\tfrac{2}{225} + \tfrac{\sqrt{6}}{75} \\
    \tfrac{2}{5} + \tfrac{\sqrt{6}}{10} &
    \tfrac{37}{225} + \tfrac{169\sqrt{6}}{1800} &
    \tfrac{11}{45} + \tfrac{7\sqrt{6}}{360} &
    -\tfrac{2}{225} - \tfrac{\sqrt{6}}{75} \\    
    1 & \tfrac{4}{9} - \tfrac{\sqrt{6}}{36} & \tfrac{4}{9} + \tfrac{\sqrt{6}}{36} & \tfrac{1}{9} \\ \hline
    & \tfrac{4}{9} - \tfrac{\sqrt{6}}{36} & \tfrac{4}{9} + \tfrac{\sqrt{6}}{36} & \tfrac{1}{9}
  \end{array}
\end{equation}
While one can also find the four-stage method explicitly, higher-order instances must be constructed numerically. 

The three-stage DIRK of Alexander~\cite{alexander1977diagonally} has the tableau
\begin{equation}
  \begin{array}{c|ccc}
    x & x & 0 & 0 \\
    \tfrac{1+x}{2} & \tfrac{1-x}{2} & x & 0 \\
    1 & y & z & x \\ \hline
      & y & z & x
  \end{array},
\end{equation}
where $x$ is defined to be the root of the polynomial $x^3 - 3 x^2 + \tfrac{3}{2} x - \tfrac{1}{6}$ that is between $\tfrac{1}{6}$ and $\tfrac{1}{2}$.  The values of $y$ and $z$ are given by
\[
\begin{split}
  y & = -\tfrac{3}{2} x^2 + 4x - \tfrac{1}{4} \\
  z & = \tfrac{3}{2} x^ - 5 x + \tfrac{5}{4}
\end{split}
\]

The DIRKs of~\cite{biswas2023design, ketcheson2020dirk} with high weak stage order are all numerically derived.  The four-stage method with formal weak stage order three is given by the tableau
\begin{equation}
  \begin{array}{c|cccc}
    0.13756544 & 0.13756544&  0.        &  0.        &  0. \\
    0.80179012 & 0.56695123&  0.23483889&  0.        &  0. \\
    2.33179673 & 1.08354073&  2.96618224&  0.44915522&  0. \\     
    1. &0.59761292&  -0.43420998& -0.05305815&  0.88965521 \\ \hline
    & 0.59761292& -0.43420998 & -0.05305815 &  0.88965521
  \end{array}
\end{equation}

\begin{acks}
  The work of S.P.M. was partially supported by an NSERC Discovery Grant.
  The work of R.C.K. was partially supported by NSF 1912653 and 2410408.
\end{acks}
\appendix
\bibliographystyle{ACM-Reference-Format}
\bibliography{references}

\end{document}